\documentclass[moor,nonblindrev]{informs3} 

\usepackage[colorinlistoftodos,bordercolor=orange,backgroundcolor=orange!20,line
color=orange,textsize=scriptsize]{todonotes}
\catcode`\@=11
\catcode`\@=12
\newcommand{\tb}[1]{#1} 
\newcommand{\jgc}[1]{#1} 

\usepackage{float}
\usepackage{subcaption}
\usepackage{caption}
\usepackage{bm}
\SingleSpacedXI 

\usepackage[colorlinks=true,breaklinks=true,bookmarks=true,urlcolor=black,
     citecolor=black,linkcolor=black,bookmarksopen=false,draft=false]{hyperref}
\usepackage{nicefrac}

\usepackage{natbib}
 \bibpunct[, ]{(}{)}{,}{a}{}{,}%
 \def\bibfont{\small}%

\newcommand{\U}{\mathcal{U}}

\newcommand{\Kexp}{\mathcal{K}_{\sf exp}}

\newcommand{\cQ}{\mathcal{Q}}

\newcommand{\cI}{\mathcal{I}}
\newcommand{\cL}{\mathcal{L}}
\newcommand{\R}{\mathbb{R}}
\newcommand{\E}{\mathbb{E}}
\newcommand{\mcE}{\mathcal{E}}
\newcommand{\PP}{\mathbb{P}}
\newcommand{\X}{\mathcal{S}}
\newcommand{\A}{\mathcal{A}}
\newcommand{\N}{\mathbb{N}}

\newcommand{\tr}{^{\top}}

\DeclareMathOperator{\kl}{KL}
\DeclareMathOperator{\dom}{dom}
\DeclareMathOperator{\stc}{s.t.}
\newcommand{\opt}{^{\star}}

\TheoremsNumberedBySection  

\EquationsNumberedBySection 

\begin{document}
\RUNAUTHOR{Grand-Cl{\'e}ment and Petrik}

\RUNTITLE{On the convex formulations of Robust MDPs}

\TITLE{On the convex formulations of robust Markov decision processes}

\ARTICLEAUTHORS{%
\AUTHOR{Julien Grand-Cl{\'e}ment}
\AFF{Information Systems and Operations Management Department, HEC Paris,  \EMAIL{grand-clement@hec.fr}}
\AUTHOR{Marek Petrik}
\AFF{Department of Computer Science, University of New Hampshire,  \EMAIL{mpetrik@cs.unh.edu}} 
} 

\ABSTRACT{%
Robust Markov decision processes (MDPs) are used for applications of dynamic optimization in uncertain environments and have been studied extensively. Many of the main properties and algorithms of MDPs, such as value iteration and policy iteration, extend directly to RMDPs. Surprisingly, there is no known analog of the MDP convex optimization formulation for solving RMDPs. This work describes the first convex optimization formulation of RMDPs under the classical sa-rectangularity and s-rectangularity assumptions. By using entropic regularization and exponential change of variables, we derive a convex formulation with a number of variables and constraints polynomial in the number of states and actions, but with large coefficients in the constraints. \tb{We further simplify the formulation for RMDPs with polyhedral, ellipsoidal, or entropy-based uncertainty sets, showing that, in these cases, RMDPs can be reformulated as conic programs based on exponential cones, quadratic cones, and non-negative orthants}. Our work opens a new research direction for RMDPs and can serve as a first step toward obtaining a tractable convex formulation of RMDPs. 
}

\KEYWORDS{Markov decision processes, robust optimization, conic optimization}
\maketitle
\section{Introduction}
Markov Decision Processes~(MDPs) represent a popular approach to sequential decision-making~\citep{puterman2014markov}. In an MDP, the decision maker aims to compute a policy that maximizes the sum of rewards accumulated over a given horizon (possibly infinite). Thanks to their modeling power and tractability, MDPs have found widespread use in domains that range from reinforcement learning~\citep{sutton2018reinforcement}, to electricity bidding~\citep{song2000optimal}, to simulating clinical decision-making~\citep{bennett2013artificial}, and to managing inventories~\citep{porteus1990stochastic}.

The majority of algorithms for solving MDPs are based on one of three main methods: 1) \emph{Value Iteration}~(VI), 2) \emph{Policy Iteration}~(PI), and 3) \emph{Linear Programming}~(LP)~\citep{puterman2014markov}. These methods were developed independently. VI traces its origin to the seminal work of \cite{bellman1966dynamic}, PI was developed by \cite{howard1960dynamic}, and the LP formulation of MDPs dates as far back as \cite{d1960probleme}. Despite having different origins, these algorithms share deep connections. Importantly, they all rely on the contraction and monotonicity properties of the dynamic programming Bellman operator associated with the MDP. Several additional connections have been established. PI can be interpreted as an implementation of the x algorithm with a block-wise pivot rule applied to the LP formulation of MDPs~\citep{ye2011simplex,ye2005new}. VI can be interpreted as a version of PI that approximates the value functions of incumbent policies using a single computation of the Bellman operator, a special case of modified PI~\citep{scherrer2015approximate}. Finally, VI and PI can be seen as gradient descent and Newton's method, respectively, applied to the residual of the Bellman operator~\citep{Filar1996,Goyal2021,grand2021convex}.

MDPs require that important model parameters, such as transition probabilities or rewards, be estimated from empirical data, which may lead to errors in the nominal rewards and transition probabilities. It is well-documented that ignoring such statistical errors can lead to severe deteriorations in the performances of the policies, both for synthetic instances~\citep{nilim2005robust,delage2010percentile} and real problems. For instance, in healthcare applications of MDPs, it may be hard to estimate the exact values of the patients' dynamics~\citep{zhang2017robust,goh2018data,grand2020robust}.
Robust MDPs~(RMDPs) ameliorate the effects of parameter errors by considering a pessimistic formulation that computes a robust policy. An optimal robust policy seeks to maximize the worst-case return over a set of plausible values of the uncertain parameters, called the {\em uncertainty set}.

For an RMDP to be tractable, its uncertainty set must satisfy certain structural assumptions. The most common such assumptions are sa-rectangularity~\citep{iyengar2005robust}, s-rectangularity~\citep{wiesemann2013robust}, k-rectangularity~\citep{mannor2016robust}, and r-rectangularity~\citep{goyal2022robust}. Interestingly, most structural properties carry over from MDPs to RMDPs with rectangular uncertainty sets. In particular, VI extends to \emph{robust value iteration}~\citep{nilim2005robust,iyengar2005robust,wiesemann2013robust} and PI extends to {\em robust policy iteration}~\citep{hansen2013strategy,ho2021partial}. Other important MDP properties that extend to RMDPs include the optimality of stationary policies~\citep{wiesemann2013robust}, Blackwell optimality, and Pontryagin's maximum property~\citep{feinberg2012handbook,goyal2022robust}. 

Despite the significant attention that RMDPs have attracted in recent years, there is no known convex optimization formulation for RMDPs, analogous to the linear programming formulation of MDPs. This is surprising because most favorable properties of MDPs extend to RMDPs, and robust variants of value iteration and policy iteration are direct extensions of their non-robust counterparts. Combined with the fact that RMDPs can be solved in polynomial time (for a fixed discount factor), it would be natural to believe that a convex formulation of the RMDP problem exists. \tb{In fact, this research question is mentioned in several seminal papers, but previous attempts to generalize the linear program formulation of MDPs to RMDPs have failed, e.g. equation (3.2) in~\cite{iyengar2005robust} and section 2.1 and in \cite{condon1990algorithms}.} One could, of course, first compute the optimal value function using value iteration and then use it to formulate the optimization problem. Such a contrived approach, however, neither offers insights into RMDP properties nor can be used to derive more efficient algorithms. 

Indeed, there are several reasons why deriving convex optimization formulations of RMDPs is an important research question. First, as in MDPs, the convex formulation can provide crucial insights into the structure of optimal policies and value functions. For instance, the state-action occupancy frequencies naturally appear as dual variables~\citep{puterman2014markov}, and the optimization formulations are a convenient tool for sensitivity analysis. Second, convex optimization formulations enable efficient solution algorithms. For instance, LP formulation of MDPs can be solved efficiently using interior-point methods~\citep{ben2001lectures} and can benefit from any advances in these general optimization algorithms~\citep{lee2014path,cohen2021solving}. These general-purpose interior point solvers often outperform value and policy iteration. They are also the only known polynomial-time algorithm for solving MDPs when the discount factor is not fixed.

\tb{Closely related to RMDPs, the problem of finding a convex formulation for perfect information mean-payoff stochastic games has been studied extensively~\citep{neyman2003stochastic,solan2015stochastic}. 
It is, in fact, one of the major open questions in this field, as there are no known polynomial-time algorithms for solving these types of games~\citep{zwick1996complexity}. Several optimization formulations have been proposed for SGs, including a non-convex quadratic program~\citep{condon1990algorithms}, a linear program for parity games~\citep{schewe2009parity}, semi-definite programming formulations~\citep{allamigeon2018solving} or convex formulations~\citep{boros2017convex}, the last three with exponentially large coefficients in the constraints.}

In this paper, we derive the first convex formulations that can approximate discounted infinite-horizon RMDP solutions to an arbitrary precision. \tb{We argue that the main difficulty for obtaining a convex formulation of RMDPs is the lack of convexity of the robust Bellman operator, coming from its saddle-point expression, and we overcome this difficulty using an appropriate regularization and a change of variables.} We focus on sa-rectangular RMDPs~\citep{iyengar2005robust} and s-rectangular RMDPs~\citep{wiesemann2013robust}. Our {\bf main contributions} are as follows.

\paragraph{Convex formulation for RMDPs.} We construct a convex formulation for computing an optimal policy for rectangular RMDPs. To the best of our knowledge, this work is the first successful attempt to obtain such a formulation. Our reformulation is based on an entropic regularization of the robust Bellman operator, combined with an exponential change of variables. The solution of the resulting convex program provides an approximate solution to the RMDP problem, which converges to the optimal RMDP solution as the entropic regularizer scales to 0. We also derive a simplified convex formulation for RMDPs with uncertainty sets described by finitely many convex constraints. \tb{For the most common models of uncertainty, such as polyhedral, ellipsoidal and relative entropy-based uncertainty, we show that RMDPs can be formulated as {\em conic} optimization programs over exponential cones, quadratic cones and non-negative orthants}. 

\paragraph{New perspectives on regularized operators.} The idea of adding regularization to dynamic programming has been used extensively in the field of reinforcement learning~\citep{neu2017unified,geist2019theory} \tb{and recently for RMDPs~\citep{kumar2022efficient,Derman2021}.}
While regularization has been used to develop variants of value iteration and policy iteration and to provide new understandings of robustness in parameter deviations, our results provide a completely new perspective on regularization. We show in this work that, surprisingly, adding regularization to RMDPs makes it possible to formulate their objectives as convex optimization problems after an appropriate change of variables. Therefore, our results provide the first evidence that regularization can be used to ensure the {\em hidden convexity} of a potentially non-convex optimization problem. 

\paragraph{A new research direction for RMDPs.} Finally, we emphasize that one of the crucial contributions of this paper is to shed light on the problem of finding a convex formulation of RMDPs. Robust value iteration and robust policy iteration have been known for two decades, and an important part of the recent literature on RMDPs has focused on obtaining efficient implementations of variants of these algorithms~\citep{kaufman2013robust,ho2021partial,Derman2021}. While finding a convex formulation has been in the spotlight in the literature on stochastic games, this problem has not been thoroughly studied in the literature on RMDPs. We assert that this is an important open problem in the field of RMDPs too. Any algorithms for solving convex optimization programs can be combined with a convex formulation to obtain new algorithms for RMDPs, potentially leading to more efficient methods based on the most recent advances in convex optimization. Our results open new promising avenues for future research, and we intend our work to provide the first building step toward obtaining a tractable convex formulation for RMDPs.

\vspace{2mm}
{\em Outline.} The rest of the paper is organized as follows. We conclude this section with important notations and definitions. We conduct a brief literature review in Section~\ref{sec:lit-review}. We provide some background on MDPs and RMDPs in Section~\ref{sec:preliminaries}. We present our convex formulation for general sa-rectangular RMDPs in Section~\ref{sec:formulation-regularized-policy}. \tb{ In Section~\ref{sec:more-concise-reformulation}, we describe a more concise convex formulation under additional assumptions and we obtain a conic program for various popular uncertainty sets. Finally, we present our results for s-rectangular uncertainty in Section~\ref{sec:s-rectangular}.}  

\subsection{Notation.}
We denote vectors with a lowercase bold font, such as $\bm{x} \in \R^n$ and matrices with an uppercase bold font, such as $\bm{X} \in \R^{n\times  n}$. The individual components of a vector are indicated using a subscript, such as $x_s$ for $s = 1, \dots, n$.  For any finite set $\X$, we write $|\X|$ for the cardinality of $\X$ and $\Delta(\X)$ for the x over $\X$, defined as the set of probability distribution over $\X$:
\[
  \Delta(\X) = \left\{ \bm{p} \in \R^{|\X|} \; | \; p_{s} \geq 0, \forall \; s \in \X, \sum_{s \in \X} p_{s}=1\right\}.
\]
For $m \in \N$ we write $[m] = \{1,...,m\}$.
We reserve the notation $\bm{e}$ for the vector $\bm{e}=(1,...,1)$, its dimension depending of the context. We write $\R_{+}^{*}=(0,+\infty)$. For $\bm{x},\bm{y} \in \R^{\X}$, we write $\bm{x} \leq \bm{y}$ for the set of component-wise inequalities $x_{s} \leq y_{s}, \forall \; s \in \X$. An operator $T\colon \R^{\X} \rightarrow \R^{\X}$ is monotone if $\bm{x} \leq \bm{y} \Rightarrow T(\bm{x}) \leq T(\bm{y})$, and it is a contraction for the $\ell_{\infty}$ norm if there exists $\gamma \in [0,1)$ such that $\| T(\bm{x}) - T(\bm{y}) \|_{\infty} \leq \gamma \| \bm{x} - \bm{y}\|_{\infty}$.
A function $f\colon  \R^{\X} \rightarrow \R$ is said to be convex if for any $\theta \in [0,1]$, for any two $\bm{x},\bm{y} \in \R^{\X}$, we have $f(\theta \bm{x} + (1-\theta)\bm{y}) \leq \theta f(\bm{x}) + (1-\theta)f(\bm{y}).$ A set $\mcE$ is \emph{convex} if for any $\theta \in [0,1]$, for any two $\bm{x},\bm{y} \in \mcE$, we have $\theta \bm{x} + (1-\theta)\bm{y} \in \mcE.$ A constraint $f(\bm{x}) \leq \bm{0}$ is said to be a {\em convex constraint} if the set $\{ \bm{x} \in \R^{\X} \; | \; f(\bm{x}) \leq \bm{0} \}$ is convex. An optimization problem is said to be a {\em convex formulation} if it minimizes a convex function (alternatively, maximizes a concave function) over a closed, convex set. The relative interior of a set $\mcE$ is denoted as ${\sf relint}(\mcE)$. We use the convention that $0 \times  \log(0) = 0$ and $0/0=1$. We use $\stc$ for {\em such that}.
\jgc{
A function $f:\R^{n} \rightarrow \R \cup \{+\infty,-\infty\}$ is {\em proper} if $f(\bm{x})>-\infty$ for all $\bm{x} \in \R^{n}$ and $f(\bm{x})<+\infty$ for at least $\bm{x} \in \R^{n}$. A function $f$ is {\em closed} if $f$ is lower semicontinuous and either $f(\bm{x}) > -\infty$ for all $\bm{x} \in \R^{n}$ or $f(\bm{x}) =- \infty$ for all $\bm{x} \in \R^{n}$. Its domain is $\dom(f):= \{ \bm{x} \in \R^{n} \; | \; f(\bm{x}) < +\infty\}$.
}

\section{Literature review}\label{sec:lit-review}

Our results combine results from several different fields of research, which we summarize next.

\paragraph{RMDPs.} RMDPs use a max-min approach to compute policies that are immune to model errors, as is common in the broader robust optimization domain~\citep{Ben-Tal2009}. While the concept of computing the best policy for a worst plausible model has a long history in decision making and MDPs~\citep{Scarf1958,satia1973markovian,White1994,givan1997bounded}, it has only garnered more widespread attention in the last two decades. The modern incarnation of RMDPs was introduced in \cite{nilim2005robust} and \cite{iyengar2005robust} along with the assumption that their uncertainty sets are sa-rectangular. With sa-rectangular sets, the adversarial nature can choose the worst model for each state and action independently. The concept of sa-rectangularity has been further extended to s-rectangularity~\citep{wiesemann2013robust}, k-rectangularity~\citep{mannor2016robust}, and r-rectangularity~\citep{goh2018data,goyal2022robust}. Without making a rectangularity assumption, even evaluating the worst-case return of a known policy is NP-hard~\citep{wiesemann2013robust}.

Rectangular RMDPs satisfy many of the same properties as MDPs. \emph{Value iteration}, \emph{policy iteration}, and \emph{modified policy iteration}, standard MDP algorithms, have been adapted to rectangular RMDPs requiring only minor modifications~\citep{ho2022robust,kaufman2013robust,grand2021scalable}. Robust value functions can also be computed efficiently for many types of uncertainty sets, such as ones based on $\ell_{1},\ell_{2}$ or $\ell_{\infty}$ norms~\citep{iyengar2005robust,ho2018fast,behzadian2021fast,Derman2021,kumar2022efficient} or $\phi$-divergence~\citep{ho2022robust}. Surprisingly, as observed by \cite{iyengar2005robust} and others in analogous contexts~\citep{hansen2013strategy}, attempts to generalize the MDP linear program formulation to RMDPs invariably result in non-convex optimization problems. The lack of a convex formulation is all the more surprising given that the standard minimax duality holds for rectangular RMDPs~\citep{wiesemann2013robust}, hinting at an underlying convexity.

An interesting special case of RMDPs assumes that the transition function is known and the rewards are uncertain~\citep{Regan2009}. Such RMDPs arise in inverse reinforcement learning which aims to imitate an unknown policy of an expert acting in a known environment~\citep{Brown2020,Javed2021}. Because the return of a policy is linear in the rewards, such reward-uncertain RMDPs are tractable and admit convex formulations without requiring any rectangularity assumptions~\citep{Regan2009,Brown2020}.

\paragraph{Regularized MDPs.}
Regularized MDPs augment the standard Bellman operator with functions that penalize the deviation of the value function or the policy from some baseline values. Numerous regularization functions have been proposed in reinforcement learning with diverse goals, including encouraging exploration~\citep{sutton2018reinforcement,Asadi2017}, risk-aversion~\citep{borkar2002q,Howard1972,Marcus1997}, or the differentiability with respect to the policy~\citep{neu2017unified}. In some cases, RMDPs are closely related to specific regularized MDPs~\citep{Derman2021,kumar2022efficient}. \tb{In particular, \cite{Derman2021} introduce {\em twice-regularized} Bellman operators, where both the policies and the value vectors are regularized, which allows for faster evaluation of the robust Bellman operator}. \tb{The connection between regularization and uncertainty in rewards is also discussed in \cite{husain2021regularized,eysenbach2021maximum}, with a focus on entropic regularization, a setting close to ours, except that we focus on uncertainty in transition probabilities.}
All these works revolve around providing (regularized) value iteration and (modified, regularized) policy iteration style algorithms for regularized MDPs.

To the best of our knowledge, our work is the first to prove that regularization, coupled with a change of variable, is also helpful in obtaining a convex formulation of RMDPs. \tb{We also note that combining regularization and changes of variables have appeared in the literature on inverse reinforcement learning \citep{lacotte2019risk,garg2021iq}, but the initial objective function in inverse reinforcement learning is already convex-concave, and the regularization and changes of variables are used to obtain closed-form solutions of the inner optimization problems, rather than inducing convexity properties.}

\paragraph{Convex formulations of zero-sum games.}
The saddle point formulation inherent to RMDPs makes them very similar to zero-sum stochastic games~\citep{neyman2003stochastic}. Whether one can formulate a stochastic game as a convex optimization problem is a natural question that also remains to be fully answered. Objectives that include mean-payoff parity games, a  model equivalent to stochastic games, can be formulated as linear programs~\citep{schewe2009parity}. This formulation has exponentially large coefficients, as do other related formulations of stochastic games~\citep{boros2017convex}. Another approach, based on semidefinite programs over a set of real series and tropical geometry, relates a convex feasibility problem and the positivity of the value of a zero-sum mean-payoff stochastic game~\citep{allamigeon2018solving}. We note that all these results apply to {\em mean-payoff} stochastic games, and they do not apply to discounted RMDPs. Indeed, \cite{schewe2009parity,boros2017convex,allamigeon2018solving} focus on {\em mean-payoff} and on a finite number of actions, whereas we focus on discounted returns and RMDPs may involve a minimization over an infinite number of transition probabilities.

\paragraph{Hidden convexity.} Since convexity is a key property in optimization, there has been many efforts to reformulate non-convex optimization programs into convex optimization programs with the same solutions. These efforts can be traced back to \cite{dantzig1963linear}, which shows a certain form of {\em hidden convexity} for the problem of generalized linear programming. More recent efforts include \cite{ben1996hidden} for quadratic programming and \cite{ben2011hidden} for partially separable optimization programs. The authors in \cite{gorissen2022hidden} study the hidden convexity in a class of problems with bilinear constraints, a setting close to RMDPs. Unfortunately, their results do not apply in our setting because they require some ``column-wise'' properties for the constraints, which do not hold for RMDPs.

\tb{There are also recent work studying the {\em convexification} of general optimization problems under suitable conditions, but the resulting optimization problems are convex only in specific parameter regimes~\citep{dvijotham2014universal}.
Closer to our setting, \cite{zhang2020variational} study MDPs with objectives that are concave in the discounted state-action occupancy measures, thus generalizing the discounted return, and then they exploit the ``hidden convexity'' property induced by the same problems viewed as optimization problems over the set of policies. We note that in contrast to our results, the optimization problems studied in \cite{zhang2020variational} are constructed directly from a hidden convex function of the state-action occupancy measures, whereas our work uncovers a hidden convexity property for the classical robust Bellman operator operating over the space of value functions.}

\section{Preliminaries on robust sequential decision-making}\label{sec:preliminaries}
\subsection{Markov Decision Processes}\label{sec:mdp}
A Markov Decision Process (MDP) is a tuple $\left(\X,\A,\bm{r},\bm{\alpha},\gamma,\bm{P}\right)$, where the {\em finite} sets $\X$ and $\A$ represent states and actions, respectively. At every decision period $t \geq 0$, the decision maker is in a given state $s \in \X$, chooses an action $a \in \A$, obtains an instantaneous reward $r_{sa} \in \R$ and transitions to the next state $s'$ at period $t+1$ with probability $P_{sas'}$, where $\bm{P} \in \left(\Delta(\X)\right)^{\X \times \A}$ is the transition function. The discount factor is $\gamma \in (0,1)$ and the initial distribution $\bm{\alpha} \in \Delta(\mathcal{S})$ determines the initial state.  We assume that $\X$ and $\A$ are finite sets: $|\X| < + \infty, |\A| < + \infty$. We also assume that all rewards are non-negative: $r_{sa} \geq 0$ for each $(s,a) \in \X \times \A$.
This assumption makes some of our bounds more convenient and is without loss of generality because adding a constant to the rewards does not change the optimal solution.

The goal of the decision maker in our setting is to find a policy $\pi$ to maximize the discounted \emph{return} $R(\pi,\bm{P})$, which is defined as the infinite-horizon discounted expected sum of rewards: 
 \begin{equation*}
 R(\pi,\bm{P})=\E^{\pi, \bm{P}} \left[ \sum_{t=0}^{\infty} \gamma^{t}r_{S_{t}A_{t}} \; \bigg| \; S_{0} \sim \bm{\alpha} \right],
 \end{equation*}
 where $(S_{t},A_{t})$ is the pair of state-action random variables visited at time $t \in \N.$
An optimal policy can be chosen \emph{stationary}~\citep{puterman2014markov}, i.e., it can be represented as a map $\pi \colon \X \to \Delta(\A)$ where $\Delta(\A)$ is the x over the set $\A$, and the action at state $s$ is always chosen following the same probability distribution $\bm{\pi}_{s} \in \Delta(A)$. There also always exists an optimal \emph{deterministic} policy which assigns the entire probability mass to a single action. We write $\Pi$ for the set of all stationary policies: $\Pi = \left(\Delta(\A)\right)^{\X}$. With these notations, the goal of the decision-maker is to find a policy $\pi\opt$ that solves the following optimization program:
\begin{equation}\label{eq:MDP}
    \max_{\pi \in \Pi} R(\pi,\bm{P}).
\end{equation}

For $\epsilon >0$, a policy $\pi \in \Pi$ is called \emph{$\epsilon$-optimal} if $R(\pi\opt,\bm{P}) - \epsilon \leq R(\pi,\bm{P}) \leq R(\pi\opt,\bm{P}).$ Note that~\eqref{eq:MDP} is not a convex optimization problem in the variable $\pi$. However, an optimal policy can be computed efficiently with value iteration, with policy iteration, by solving a linear program~\citep{puterman2014markov}, and even by gradient descent~\citep{Bhandari2021}. Although these approaches have been developed mostly independently~\citep{d1960probleme,howard1960dynamic,bellman1966dynamic}, they are based on computing the fixed point of the dynamic programming operator associated with the MDP, called the Bellman operator.

\paragraph{The Bellman operator.}
Let us consider some fixed transition probabilities $\bm{P} \in \left(\Delta(\X)\right)^{\X \times \A}$. For each policy $\pi \in \Pi$, we define $\bm{v}^{\pi} \in \R^{\X}$ the {\em value function}:
\begin{equation} \label{eq:value-policy-nominal}
     v_{s}^{\pi} \;=\; \E^{\pi, \bm{P}} \left[ \sum_{t=0}^{\infty} \gamma^{t}r_{S_{t}A_{t}} \; \bigg| \; S_{0} =s \right].
\end{equation}
 Algebraic manipulation shows that $R(\pi,\bm{P}) = \bm{\alpha}^{\top}\bm{v}^{\pi}$~\citep{puterman2014markov}. 
 Define the \emph{Bellman operator} $T_{\bm{P}}\colon \R^{\X} \rightarrow \R^{\X}$ for the nominal MDP with fixed transition probabilities $\bm{P}$ as follows:
\begin{equation}\label{eq:Bellman-operator-nominal}
    T_{\bm{P}}(\bm{v})_{s} = \max_{a \in \A} \; r_{sa} + \gamma \bm{P}_{sa}^{\top}\bm{v}, \quad  \forall s \in \X,\, \forall \; \bm{v} \in \R^{\X}.
\end{equation}
Similarly, we can define for any Markov policy $\pi$ a \emph{Bellman evaluation operator} $T_{\bm{P}}^{\pi}\colon \R^{\X} \to \R^{\X}$ such that 
\begin{equation}\label{eq:Bellman-operator-nominal-policy}
    T_{\bm{P}}^{\pi}(\bm{v})_{s} = \sum_{a \in \A} \pi_{sa} \cdot \left(r_{sa} + \gamma \bm{P}_{sa}^{\top}\bm{v}\right), \quad  \forall s \in \X,\, \forall \bm{v} \in \R^{\X}.
\end{equation}
A thorough exposition of the properties of MDPs and of the classical methods to solve MDPs~(VI, PI, LP formulation) can be found in section~6 in \cite{puterman2014markov}. In particular, as shown in section~6.3 in \cite{puterman2014markov}, the operator $T_{\bm{P}}$ is a contraction for the $\ell_{\infty}$ norm:
    \[\| T_{\bm{P}}(\bm{v}) - T_{\bm{P}}(\bm{v}') \|_{\infty} \leq \gamma \cdot \| \bm{v} - \bm{v}'\|_{\infty}, \quad \forall \; \bm{v},\bm{v}' \in \R^{\X},\]
and $T_{\bm{P}}$ is monotone:
\[ v_{s}\leq v_{s}', \forall s \in \X \quad \Longrightarrow\quad  T_{\bm{P}}(\bm{v})_{s} \leq T_{\bm{P}}(\bm{v}')_{s},\forall s \in \X.\]
Because $T_{\bm{P}}$ is a contraction, it has a unique fixed point.
The following proposition summarized important properties of the fixed-point of $T_{\bm{P}}$.
\begin{proposition}\label{prop:Bellman-equation-MDP}
Let $\bm{v}\opt$ be the unique fixed-point of the Bellman operator $T_{\bm{P}}$. Then $\bm{v}\opt$ is the value function of an optimal policy $\pi\opt$ in the MDP problem~\eqref{eq:MDP}. Moreover, $\pi\opt$ can be chosen deterministic and attains the maximum on the right-hand side of $T_{\bm{P}}(\bm{v}\opt)_{s}$ as in Equation~\eqref{eq:Bellman-operator-nominal} for each state $s\in \X$.
\end{proposition}
Proposition~\ref{prop:Bellman-equation-MDP} can be traced back to the seminal works of Shapley~\citep{shapley1953stochastic} and Bellman~\citep{bellman1966dynamic}.
The value function $\bm{v}^{\pi}$ defined in~\eqref{eq:value-policy-nominal} for a policy $\pi$ satisfies a property analogous to Proposition~\ref{prop:Bellman-equation-MDP} and is a fixed point of the operator $T_{\bm{P}}^{\pi}$, which is also a contraction for the $\ell_{\infty}$ norm. 
\tb{Both Value Iteration (VI) and Policy Iteration (PI) are based on Proposition~\ref{prop:Bellman-equation-MDP}, and efficiently solve the MDP problem~\eqref{eq:MDP}. We refer to Appendix~\ref{app:alg-mdps} for more details on algorithms for solving nominal MDPs.}
\paragraph{Contraction lemma and linear programming.}
To introduce the linear programming formulation of MDPs, we start with the following {\em contraction lemma}, which plays a critical role in all the results in this paper.
\begin{lemma}[Contraction lemma]\label{lem:contraction-program}
Let $F\colon  \R^{\X} \rightarrow \R^{\X}$ be a monotone contraction operator, $g\colon \R^{\X} \rightarrow \R$ be a component-wise non-decreasing function, and $\bm{v}\opt$ the unique fixed-point of $F$.
Then
\begin{equation}\label{eq:optimization-fixed}
  g(\bm{v}\opt) \;=\;  \min \{ g(\bm{v}) \; | \; \bm{v} \geq F(\bm{v}) \} \; =\max \{ g(\bm{v}) \; | \; \bm{v} \leq F(\bm{v}) \} \; .
\end{equation}
Moreover, if $g$ is component-wise increasing then $\bm{v}\opt$ is the unique solution to both optimization problems in~\eqref{eq:optimization-fixed}.
\end{lemma}
Versions of Lemma~\ref{lem:contraction-program} are presented in \cite{nilim2005robust} and \cite{goyal2022robust}. For the sake of completeness, we provide its proof in Appendix~\ref{app:proof-contraction-lemma}. Lemma~\ref{lem:contraction-program} provides a powerful tool to construct convex formulations. In particular, if $F_{s}$ is convex for each $s \in \X$ and $g$ is convex, then $\min \{ g(\bm{v}) \; | \; \bm{v} \geq F(\bm{v}) \}$ is a convex program. Similarly, if $F_{s}$ is concave for each $s \in \X$ and $g$ is concave, then $\max \{ g(\bm{v}) \; | \; \bm{v} \leq F(\bm{v}) \}$ is a convex program. We will use Lemma~\ref{lem:contraction-program} as a building block to construct our convex formulation of RMDPs.

In particular, applying Lemma~\ref{lem:contraction-program} with $F(\bm{v})=T_{\bm{P}}(\bm{v})$ and $g(\bm{v}) = \bm{\alpha}\tr \bm{v}$ yields that $\bm{v}\opt$ is an optimal solution to $\min \{ \bm{\alpha}\tr\bm{v} \; | \; \bm{v} \geq T_{\bm{P}}(\bm{v})\}.$ This directly gives a convex optimization formulation of MDPs, since for each component $s \in \X$, the map $\bm{v} \mapsto T(\bm{v})_{s}$ is convex as the maximum of linear forms (section~3.2.3, \cite{boyd2004convex}).
Additionally, note that 
\[
  \bm{v} \geq T_{\bm{P}}(\bm{v})
  \quad \iff \quad
  v_{s} \geq \max_{a \in \A} r_{sa} + \gamma \bm{P}_{sa}^{\top}\bm{v}, \forall \; s \in \X
  \quad \iff\quad
  v_{s} \geq r_{sa} + \gamma \bm{P}_{sa}^{\top}\bm{v}, \forall \; (s,a) \in \X \times \A.
\]
We conclude that $\bm{v}\opt$ is a solution to the following linear program with $|\X|$ variables and $|\X| \times |\A|$ constraints:
\begin{equation}\label{eq:linear-program-mdp-primal}
    \begin{aligned}
        \min_{\bm{v}}\quad  & \bm{\alpha}\tr\bm{v}\\
        \stc \quad &  v_{s} \geq r_{sa} + \gamma \bm{P}_{sa}^{\top}\bm{v}, \quad  \forall \; (s,a) \in \X \times \A, \\
        & \bm{v} \in \R^{\X}.
    \end{aligned}
\end{equation}
Therefore, while MDPs are not convex programs in the formulation~\eqref{eq:MDP} where the objective is $\pi \mapsto R(\pi,\bm{P})$ and the variable is the policy $\pi \in \Pi$, the formulation~\eqref{eq:linear-program-mdp-primal} provides a linear program with variables $\bm{v} \in \R^{\X}$, \tb{whose dual variables are interpreted as state-action occupancy frequency (see Appendix~\ref{app:alg-mdps}). An optimal policy $\pi\opt$ can directly be recovered from an optimal solution $\bm{v}\opt$ of~\eqref{eq:linear-program-mdp-primal}, as the policy $\pi\opt$ such that $T_{\bm{P}}(\bm{v}\opt) = T^{\pi\opt}_{\bm{P}}(\bm{v}\opt)$.} 
\subsection{Robust Markov Decision Processes}\label{sec:rmdp}
Robust MDPs~(RMDPs) generalize MDPs to allow for uncertain transition function $\bm{P}$.
Specifically, RMDPs replace a known $\bm{P}$ by an {\em uncertainty set} $\U \subseteq \left(\Delta(\X)\right)^{\X \times \A}$ of possible values for the uncertain transition function $\bm{P}$. The decision maker's goal in an RMDP is to solve the following saddle point problem:
\begin{equation}\label{eq:robust-mdp}
  \max_{\pi \in \Pi} \; \min_{\bm{P} \in \U} R(\pi,\bm{P}).
\end{equation}

An RMDP problem~\eqref{eq:robust-mdp} is only tractable under some conditions on the set $\U$. A common condition is sa-rectangularity, which assumes that the marginals $\bm{P}_{sa}$ can be chosen independently across each pair $(s,a) \in \X \times \A$~\citep{iyengar2005robust,nilim2005robust}. We make use of this structural assumption in Section~\ref{sec:formulation-regularized-policy} and Section~\ref{sec:more-concise-reformulation}, \tb{before extending our results to s-rectangular RMDPs~\citep{wiesemann2013robust} in Section~\ref{sec:s-rectangular}.}
\begin{assumption}
The uncertainty set $\U$ is sa-rectangular:
\begin{equation*}
    \U = \times_{(s,a) \in \X \times \A} \, \U_{sa}, \quad  \U_{sa} \subseteq \Delta(\X). 
\end{equation*}
\end{assumption}
Several other conditions on $\mathcal{U}$ that make~\eqref{eq:robust-mdp} tractable have been studied. These other conditions include $s$-rectangularity~\citep{wiesemann2013robust} and $r$-rectangularity~\citep{goyal2022robust}. We will later focus on s-rectangular RMDPs in Section~\ref{sec:s-rectangular}.
For any sa-rectangular RMDP one can generalize the Bellman operator in~\eqref{eq:Bellman-operator-nominal} to the {\em robust Bellman operator} $T\colon \R^{S} \rightarrow \R^{S}$ as 
\begin{equation}\label{eq:bellman-operator-robust}
  T(\bm{v})_{s} =\max_{a \in \A} \; \min_{\bm{p} \in \U_{sa}} \{ r_{sa} + \gamma \cdot  \bm{p}\tr\bm{v}\},
  \quad \forall \; s \in \X, \forall \; \bm{v} \in \R^{\X}.
\end{equation}

The robust operator $T$ is a monotone contraction for the $\ell_{\infty}$ norm~\citep{iyengar2005robust,nilim2005robust}. Therefore, we can extend Proposition~\ref{prop:Bellman-equation-MDP} to RMDPs as follows.
\begin{proposition}\label{prop:Bellman-equation-RMDP}
Let $\bm{v}\opt$ be the unique fixed-point of the robust Bellman operator $T$. Then $\bm{v}\opt$ is the value function of an optimal policy $\pi\opt$ in the RMDP problem~\eqref{eq:robust-mdp}. Moreover, $\pi\opt$ can be chosen deterministic and attains the maximum on the right-hand side of $T(\bm{v}\opt)_{s}$ as in Equation~\eqref{eq:bellman-operator-robust} for each state $s\in \X$.
\end{proposition}
We can also define for any stationary policy $\pi$ a \emph{robust Bellman evaluation operator} $T^{\pi}\colon \R^{\X} \to \R^{\X}$ such that 
\begin{equation}\label{eq:Bellman-operator-robust-policy}
    T^{\pi}(\bm{v})_{s} = \sum_{a \in \A} \pi_{sa} \cdot \min_{p\in \U_{sa}} \left(r_{sa} + \gamma \bm{p}^{\top}\bm{v}\right), \quad  \forall s \in \X,\, \forall \; \bm{v} \in \R^{\X}.
\end{equation}
The fixed point of $T^{\pi}$ is the robust value function $\bm{v}_{\U}^{\pi} \in \R^{\X}$ of a policy $\pi$, which is defined as
\begin{equation}\label{eq:value-function-robust}
     v_{\U,s}^{\pi}=\E^{\pi, \bm{P}(\pi)} \left[ \sum_{t=0}^{\infty} \gamma^{t}r_{S_{t}A_{t}} \; \bigg| \; S_{0} =s \right],  \quad \bm{P}(\pi) \in \arg \min_{\bm{P} \in \U} R(\pi,\bm{P}).
\end{equation}
The crucial difference between $T^{\pi}$ and the regular Bellman evaluation operator $T_{\bm{P}}^{\pi}$ is that $T^{\pi}$ is non-linear. Whereas one can compute the fixed point of $T_{\bm{P}}^{\pi}$ by simply solving a system of linear equations, that is impossible for $T^{\pi}$.

Proposition~\ref{prop:Bellman-equation-RMDP} plays an important role in computing an optimal policy $\pi\opt$ in an RMDP. Most methods first compute $\bm{v}\opt$ which is the unique fixed-point to the robust Bellman operator $T$. Additionally, both value iteration and policy iteration can be extended to {\em robust value iteration}~\citep{iyengar2005robust,nilim2005robust} and {\em robust policy iteration}~\citep{iyengar2005robust,hansen2013strategy,kaufman2013robust,ho2021partial} for solving RMDPs. Robust value iteration simply replaces the Bellman operator $T_{\bm{P}}$ by the robust Bellman operator $T$:
\begin{equation*}
    \bm{v}_{0} \in \R^{\X}, \quad  \bm{v}_{k+1} = T(\bm{v}_{k}), \quad \forall \; k \in \N,
\end{equation*}
while robust policy iteration replaces the sequence of value functions with a sequence of robust value functions: $\bm{v}_{0} \in \R^{\X}$, and for $k \in \N$,
\begin{equation*}
\pi_{k} \text{ such that }  T^{\pi_k}(\bm{v}_{k}) = T(\bm{v}_{k}),  
  \text{ and }
\bm{v}_{k+1}  = \bm{v}_{\U}^{\pi_{k}}.
\end{equation*}
For the sake of brevity, we omit the details on the convergence rates of robust VI and robust PI to an optimal value function. We simply note that at each iteration $k$, the robust VI now must solve a saddle point problem to evaluate $T(\bm{v}_{k})$, while robust PI requires to compute $\bm{v}^{\pi_{k}}_{\U}$, which can be computationally challenging~\citep{ho2021partial}.

To obtain a convex formulation of RMDPs, it is natural to take the same approach as for MDPs and start from the contraction lemma. In particular, for RMDPs, it still holds that
\begin{equation}\label{eq:v-star-arg-min-arg-max}
  \bm{v}\opt \in
  \arg \min \; \{ \bm{\alpha}^{\top}\bm{v} \; \vert \; \bm{v} \geq T(\bm{v}) \}
  =
  \arg \max \; \{ \bm{\alpha}^{\top}\bm{v} \; \vert \; \bm{v} \leq T(\bm{v}) \}.
\end{equation}
However, for each $s \in \X$, the map $\bm{v} \mapsto T(\bm{v})_{s}$ may not be convex nor concave, because of the saddle point formulation for the robust Bellman operator $T$ as in~\eqref{eq:bellman-operator-robust}. Because of the lack of convexity, applying Lemma~\ref{lem:contraction-program} to the robust Bellman operator $T$ does not directly lead to a convex formulation of the RMDP problem~\eqref{eq:robust-mdp}.  We present a  example of this fact below.
\begin{example}\label{ex:non-convex-T}
We consider the following instance: $\X=\{1,2\},\A=\{1,2,3\},\gamma = 0.8$, $
\bm{r}_{1} = (2,11,10),\bm{r}_{2} = (1,1,1)$.
The nominal transition probabilities are $\hat{\bm{P}}$ with
\[
\begin{array}{rlrlrl}
  \hat{\bm{P}}_{1,1}& =(0.1,0.9),
  &\hat{\bm{P}}_{1,2}&=(0.25,0.75),
  &\hat{\bm{P}}_{1,3}&=(0.4,0.6),\\
  \hat{\bm{P}}_{2,1}& =(0.5,0.5),
  &\hat{\bm{P}}_{2,2}&=(0.5,0.5),
  &\hat{\bm{P}}_{2,3}&=(0.5,0.5),
\end{array}
\]
and the uncertainty sets $\U_{sa}$ are $\ell_{\infty}$ balls centered around $\hat{\bm{P}}$:
\begin{align*}
    \U_{sa} = \{ \bm{p} \in \Delta(\X) \; | \; 0.95 \cdot \hat{\bm{p}}_{sa} \leq \bm{p} \leq 1.05 \cdot \hat{\bm{p}}_{sa}\}, \forall \; (s,a) \in \X \times \A.
\end{align*}
We take $\bm{v}_{1} = (10.5,0.85),\bm{v}_{2} = (0.5,4)$ and we write $\bm{v}_{\theta} = \theta\bm{v}_{1}+(1-\theta)\bm{v}_{2}$ with $\theta \in [0,1]$. In Figure~\ref{fig:non-convex-T}, we plot the map $\theta \mapsto T(\bm{v}_{\theta})_{1}$, which is neither concave nor convex for $s=1$. This map is piece-wise affine, with the first change of slopes corresponding to a change in (local) optimal probabilities (in $\min_{\bm{p} \in \U_{sa}}$), and the second change of slopes corresponding to a change in (local) optimal actions (in $\max_{a \in \A}$).
\begin{figure}[hbt]
\begin{center}
         \includegraphics[width=0.4\linewidth]{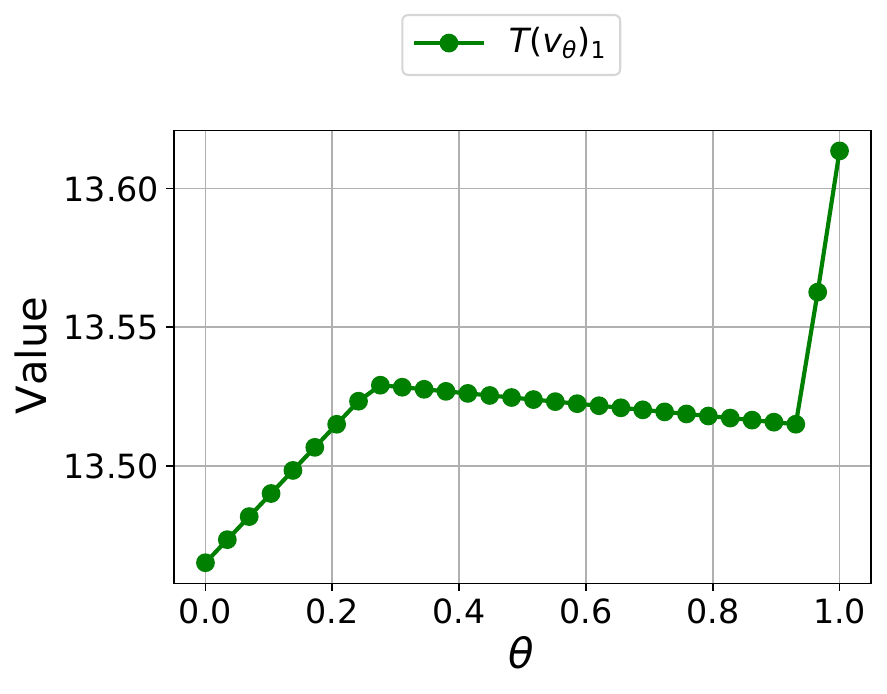}
\end{center}
  \caption{Non-concavity and non-convexity of $\theta \mapsto T_{1}(\bm{v}_{\theta})$.}
  \label{fig:non-convex-T}
\end{figure}
\end{example}
\begin{remark}[Optimistic MDPs]
In the case of {\em optimistic}  MDPs~\citep{Auer2010}, the goal is to compute the best-case possible return of a policy:
\begin{equation}\label{eq:robust-optimistic-MDP}
    \max_{\pi \in \Pi} \max_{\bm{P} \in \U} R(\pi,\bm{P}).
\end{equation}
In this case, the robust Bellman operator becomes
\[
    T^{\sf opt}(\bm{v})_{s} =\max_{a \in \A} \; \max_{\bm{p} \in \U_{sa}} \{ r_{sa} + \gamma \cdot  \bm{p}\tr\bm{v}\}, \quad \forall \; s \in \X, \forall \; \bm{v} \in \R^{\X},\]
    and $T^{\sf opt}_{s}$ is now a convex function for each $s \in \X$, as the maximum of some convex functions. Therefore, in this case a convex formulation of~\eqref{eq:robust-optimistic-MDP} is given by $\max \; \{ \bm{\alpha}\tr\bm{v} \; | \; \bm{v} \geq T^{\sf opt}(\bm{v})\}$.
\end{remark}
\section{Convex formulation based on entropic regularization}
  \label{sec:formulation-regularized-policy}
In this section, we describe a convex formulation of the RMDP problem. \tb{Our formulation is based on Lemma~\ref{lem:contraction-program}, and overcomes the non-convexity of the robust Bellman operator with regularization and changes of variables}. For the sake of simplicity, we present our results for sa-rectangular uncertainty sets, even though the results in this section extend to the case of s-rectangular uncertainty sets, as described in Section~\ref{sec:s-rectangular}. We first introduce regularized operators in Section~\ref{sec:regularized-operators} before constructing our convex formulation in Section~\ref{sec:convex-formulation}. 
\subsection{Regularized operators}\label{sec:regularized-operators}
Recall that the robust Bellman operator $T\colon  \R^{\X} \rightarrow \R^{\X}$ is defined as
\[T(\bm{v})_{s} = \max_{a \in \A} \left(r_{sa} + \gamma \min_{\bm{p} \in \U_{sa}} \bm{p}^\top \bm{v}\right) = \max_{\bm{\pi}_{s} \in \Delta(\A)}  \sum_{a \in \A} \pi_{sa}\left(r_{sa} + \gamma \min_{\bm{p} \in \U_{sa}} \bm{p}^\top \bm{v}\right) , \; \forall \; s \in \X.\]
As explained in the previous section, for each $s \in \X$, the map $\bm{v} \mapsto T(\bm{v})_{s}$ may be neither convex nor concave in $\bm{v}$, because of the saddle point formulation in the robust Bellman operator. Our convex formulation of RMDPs combines {\em entropic regularization} with an exponential change of variables. In particular, we start by introducing a \emph{regularized robust Bellman operator} $\tilde{T}\colon \R^{\X} \rightarrow \R^{\X}$, defined as
\begin{equation}\label{eq:T-tilde-definition}
    \tilde{T}(\bm{v})_{s} = \max_{\bm{\pi}_{s} \in \Delta(\A)}  \sum_{a \in \A} \pi_{sa}\left(r_{sa} + \gamma \min_{\bm{p} \in \U_{sa}} \bm{p}^\top \bm{v}\right) - \frac{1}{b} \cdot \kl(\bm{\pi}_{s},\bm{\nu}_{s}), \; \forall \; s \in \X,
  \end{equation}
  where $\nu \in \Pi$ is a policy and $\kl\colon \Delta(\mathcal{A}) \times  \Delta(\mathcal{S}) \to \mathbb{R}$ is the Kullback-Leibler divergence \tb{(sometimes called the {\em relative entropy}}):
\[
   \kl(\bm{\pi}_{s},\bm{\nu}_{s}) = \begin{cases}
        \sum_{a \in \A} \pi_{sa} \log\left(\frac{\pi_{sa}}{\nu_{sa}}\right) \quad &\text{if}  \; \bm{\pi}_s \ll \bm{\nu}_s, \\
        + \infty \quad &\text{otherwise}.
   \end{cases} 
\]
 \tb{Here, $\bm{\pi}_s \ll \bm{\nu}_s$ indicates that $\bm{\pi}_s$ is absolutely continuous with respect to $\bm{\nu}_s$ and is defined as $\nu_{sa} =0 \Rightarrow \pi_{sa}=0, \forall \; a \in \A$.
The policy $\nu$ can be used to incorporate prior knowledge on an optimal policy or to bias a solution toward $\nu_{s}$. In practice, it is typical to choose a uniform policy over the set of actions, and in the rest of the paper we assume that $\nu_{sa} = 1/|\A|, \forall \; (s,a) \in \X \times \A$.}
The idea of entropic regularization has been used extensively in reinforcement learning~\citep{szepesvari2010algorithms,neu2017unified}, but it has mostly provided variants of value iteration and policy iteration replacing the robust Bellman operator $T$ with the regularized robust Bellman operator $\tilde{T}$~\citep{geist2019theory,Derman2021,kumar2022efficient}. 
An alternative view is to incorporate the $\kl$ term in the definition of $\tilde{T}$ as in~\eqref{eq:T-tilde-definition} as a modeling assumption on the instantaneous rewards. This is the point of view adopted in the model of linearly-solvable (nominal) MDPs~\citep{todorov2006linearly}. In that sense, the operator $\tilde{T}$ can be interpreted as the Bellman operator associated with a  linearly-solvable {\em robust} MDPs, or equivalently regularized RMDPs. While it could be of interest to study this new model (e.g., optimality of deterministic Markovian policies, regularizers that lead to tractable updates for value iteration and policy iteration, etc.), we focus our efforts on finding a convex formulation of RMDPs and keep the study of robust regularized MDPs for future work. 
\begin{remark}[Connection with Policy Mirror Descent]
    \tb{We note that when there is no uncertainty, i.e., when $\U$ is a singleton, the regularized robust Bellman operator reduces to the classical {\em regularized Bellman operator} from definition~1 in \cite{geist2019theory}, which is related to the Natural Policy Gradient algorithm~\citep{kakade2001natural} and its extension as Policy Mirror Descent~\citep{lan2023policy}. In contrast, when $\U$ is not a singleton, the proximal updates from Robust Policy Mirror Descent~\citep{li2022robust} involve robust Q-functions. This differs from our regularized robust Bellman operator $\tilde{T}$ which only involves an inner minimization over $\bm{p} \in \U_{sa}$. In this paper, we will only consider a {\em fixed} regularization scalar $b$, in contrast to what the changing step sizes in first-order methods for robust MDPs~\citep{lan2023policy}. This is because we are interested in convex formulations of RMDPs rather than in designing iterative algorithms.}
\end{remark}
Crucially, the operator $\tilde{T}$ is still a monotone contraction, and its unique fixed-point $\tilde{\bm{v}}\opt$ can be seen as an approximation of the unique fixed-point $\bm{v}\opt$ of the operator $T$. In particular, we have the following proposition.
\begin{proposition}\label{prop:T-T-tilde}
 \tb{When $\nu_{sa} = 1/|\A|$ for each $(s,a) \in \X \times \A$ in the regularized operator $\tilde{T}$, then:}
\begin{enumerate}
    \item The operator $\tilde{T}$ is a monotone contraction with respect to the $\ell_{\infty}$ norm.
    \item For each $s \in \X$ and $\bm{v} \in \R^{\X}$, the operator $\tilde{T}$ brackets $T$ as
      \[
        \tilde{T}(\bm{v})_{s} \leq T(\bm{v})_{s}  \leq \tilde{T}(\bm{v})_{s} + \frac{\log(|\A|)}{b}.
      \]
    \item Let $\tilde{\bm{v}}\opt \in \R^{\X}$ the unique fixed-point of $\tilde{T}$ and $\bm{v}\opt \in \R^{\X}$ the unique fixed-point of $T$. Then 
    \[ \tilde{v}\opt_{s} \leq v\opt_{s}  \leq \tilde{v}\opt_{s} + \frac{\log(|\A|)}{b \cdot (1-\gamma)}, \quad  \forall \; s \in \X.\]
\end{enumerate}
\end{proposition}
We present a detailed proof in Appendix~\ref{app:prf:prop:T-T-tilde}. Note that analogous results for the nominal Bellman operator $T_{\bm{P}}$ can be found in section~3 in \cite{geist2019theory}.
We have the following lemma, which provides a  formulation of $\tilde{T}$.
\begin{lemma}\label{lem:kl-fenchel-conjugate}
Let $\bm{q}_{0} \in \Delta(\A)$ and $\bm{y} \in \R^{\A}$. Then we have
\begin{equation}\label{eq:kl-fenchel-conjugate}
  \max_{\bm{q} \in \Delta(\A)} \bm{q}\tr\bm{y} - \frac{1}{b} \kl(\bm{q},\bm{q}_{0})
 \; = \;
  \frac{1}{b} \log \left( \sum_{a \in \A} q_{0,a}\exp\left(b \cdot y_{a} \right) \right).
\end{equation}
\end{lemma}
Lemma~\ref{lem:kl-fenchel-conjugate} is known as a Donsker-Varadhan variational formula~\citep{Dupuis1997} and has been proven independently from KKT conditions many times~\citep{nilim2005robust,iyengar2005robust}. In convex optimization, Equality~\eqref{eq:kl-fenchel-conjugate} shows that the conjugate of the negative entropy is the log-sum-exp function (example~3.25 in \cite{boyd2004convex}). Lemma~\ref{lem:kl-fenchel-conjugate} can also be interpreted as the dual representation of the entropic risk measure in the theory of convex risk measures~\citep{follmer2002convex,Ahmadi-Javid2012}.

Applying Lemma~\ref{lem:kl-fenchel-conjugate} to $\tilde{T}$ yields
\begin{equation}\label{eq:T-tilde-formulation}
    \tilde{T}(\bm{v})_{s} = \frac{1}{b} \log \left(   \sum_{a \in \A} \frac{1}{|\A|}\exp\left(b \cdot r_{sa} + b \cdot \gamma \min_{\bm{p} \in \U_{sa}} \bm{p}^\top \bm{v} \right) \right), \quad  \forall \; s \in \X.
\end{equation}

 Based on the contraction lemma (Lemma~\ref{lem:contraction-program}), we can compute $\tilde{\bm{v}}\opt$, the unique fixed-point of $\tilde{T}$, as the unique optimal solution to the following optimization program:
\begin{equation}\label{eq:contraction-lemma-t-tilde}
    \max \left\{ g(\bm{v}) \; | \; \bm{v} \leq \tilde{T}(\bm{v})\right\}
\end{equation}
for any component-wise increasing function $g\colon \R^{\X} \rightarrow \R$. 
Unfortunately, for each $s \in \X$, the map $\bm{v}\mapsto\tilde{T}(\bm{v})_{s}$ may be neither concave nor convex in $\bm{v}$, because of the minimization over $\U_{sa}$ for each pair $(s,a) \in \X \times \A.$ This is highlighted in the next example.
\begin{example}\label{ex:non-convex-tilde-T}
We consider the same RMDP instance and empirical setup as in Example~\ref{ex:non-convex-T}. \tb{As in Proposition~\ref{prop:T-T-tilde},} we set $\nu$ to be the uniform policy: $\nu_{sa} = \nicefrac{1}{3}$ for each $s\in \X$ and $a\in \A$. Figure~\ref{fig:comparison-T-T-tilde} shows the maps $\theta \mapsto T_{1}(\bm{v}_{\theta})$ and $\theta \mapsto \tilde{T}_{1}(\bm{v}_{\theta})$ for $\theta \in [0,1]$ and for $b=5$ and $b=10$. The map $\theta \mapsto \tilde{T}_{1}(\bm{v}_{\theta})$ is clearly neither concave nor convex.
\begin{figure}[hbt]
    \begin{center}
         \includegraphics[width=0.4\linewidth]{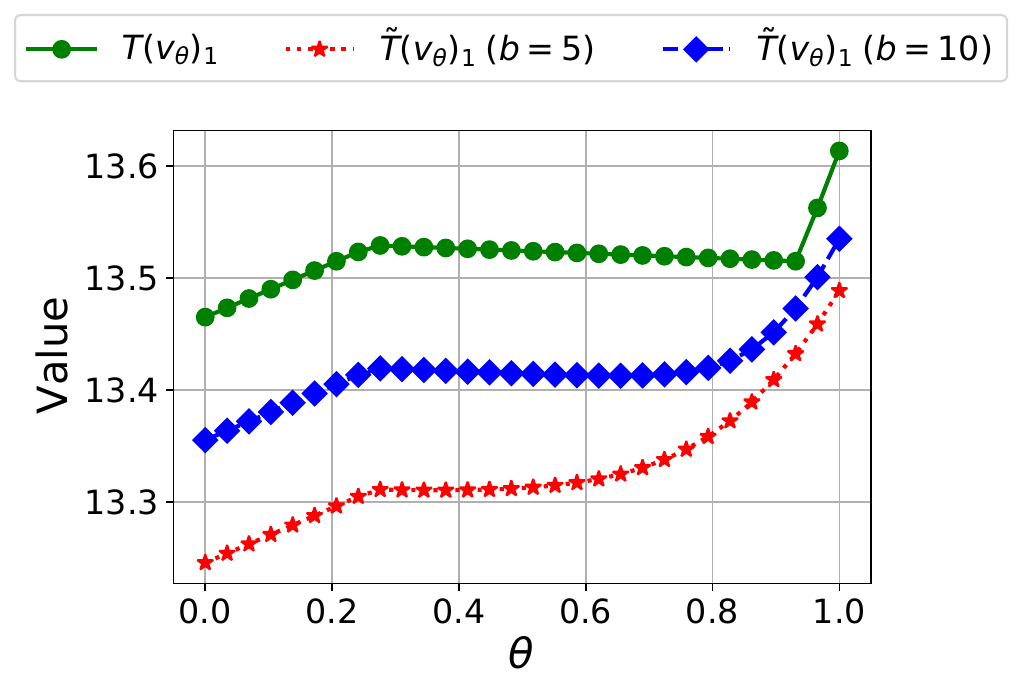}
    \end{center}
     \caption{Non-concavity and non-convexity of $\theta \mapsto T_{1}(\bm{v}_{\theta})$ and $\theta \mapsto \tilde{T}_{1}(\bm{v}_{\theta})$ for $b=5$ and $b=10$.}
    \label{fig:comparison-T-T-tilde}
  \end{figure}
\end{example}
\subsection{Exponential change of variables and convex formulation}\label{sec:convex-formulation}
We now construct our convex formulation of RMDPs, based on an exponential change of variables for the operator $\tilde{T}.$ In particular, recall that for all $\bm{v} \in \R^{\X},$
\[
  \tilde{T}(\bm{v})_{s} = \frac{1}{b} \log \left(   \sum_{a \in \A} \frac{1}{|\A|}\exp\left(b \cdot r_{sa} + b \cdot \gamma \min_{\bm{p} \in \U_{sa}} \bm{p}^\top \bm{v} \right) \right), \quad  \forall \; s \in \X.
\]
Example~\ref{ex:non-convex-tilde-T} shows that $\bm{v} \mapsto \tilde{T}(\bm{v})_{s}$ may be non-convex and non-concave for some $s \in \X,$ and therefore that the set $\{ \bm{v} \in \R^{\X} \; | \; \bm{v} \leq \tilde{T}(\bm{v})\}$ may be non-convex.
However, we will show that this set becomes convex after a  {\em exponential} change of variables. 
For $b>0$, let us write $\log_{b}$ and $\exp_{b}$ for the following ``scaled'' logarithm and exponential functions:
\[
  \log_{b}(z) = \log(z)/b, \qquad  \exp_{b}(z') = \exp(b \cdot z'), \quad  \forall \; z>0, \forall \; z' \in \R.
\]
Note that $b$ is not the base of the logarithm, it is the denominator in $z \mapsto \log(z)/b$.
We will also apply the functions $\log_{b}$ and $\exp_{b}$ component-wise: for $\bm{v} \in \R^{\X}$, we note $\exp_{b}(\bm{v})$ for the vector in $\R^{\X}$ defined as $\left(\exp(b \cdot v_{s})\right)_{s \in \X}$ and $\log_{b}=\exp_{b}^{-1}.$
With this notation, we obtain
\begin{align*}
    \tilde{T}(\bm{v})_{s} & = \log_{b}\left(  \sum_{a \in \A} \frac{1}{|\A|}\exp_{b}\left(r_{sa}\right) \exp_{b} \left(  \gamma \min_{\bm{p} \in \U_{sa}} \bm{p}^\top \bm{v} \right) \right) \\
    & =  \log_{b}\left(  \sum_{a \in \A} \frac{1}{|\A|}\exp_{b}\left(r_{sa}\right) \min_{\bm{p} \in \U_{sa}} \exp_{b} \left(  \gamma  \bm{p}^\top \bm{v} \right) \right) \\
    & =  \log_{b}\left(  \sum_{a \in \A} \frac{1}{|\A|}\exp_{b}\left(r_{sa}\right) \min_{\bm{p} \in \U_{sa}} \prod_{s' \in \X}  \left( \exp_{b} \left(v_{s'}\right)\right)^{\gamma p_{s'}} \right).
\end{align*}
In particular, we can reformulate $\tilde{T}\colon \R^{\X} \rightarrow \R^{\X}$ as
\[
  \tilde{T}(\bm{v}) = \log_{b} \,\circ\, \tilde{t} \,\circ\, \exp_{b}(\bm{v}),
\]
with the operator $\tilde{t}\colon \left(\R_{+}^{*}\right)^{\X} \rightarrow \left(\R_{+}^{*}\right)^{\X}$ defined  for each state $s \in \X$ and $\bm{x} \in \left(\R_{+}^{*}\right)^{\X}$ as
\begin{equation}\label{eq:definition-f}
  \tilde{t}(\bm{x})_{s}=  \sum_{a \in \A} \frac{1}{|\A|} \exp_{b}\left(r_{sa}\right)  \min_{\bm{p} \in \U_{sa}}  \prod_{s' \in \X} x_{s'}^{\gamma p_{s'}}.
\end{equation}
\tb{For the sake of brevity, we use $\omega_{b}(s,a) = \frac{1}{|\A|} \exp_{b}\left(r_{sa}\right)$.}
The operator $\tilde{t}$ plays a crucial role in our convex formulation of RMDPs. In particular, this operator is component-wise concave, as we show in the following lemma. Recall that $\R_{+}^{*}=(0,+\infty)$.
\begin{lemma}\label{lem:concavity-of-f}
For each $s \in \X$, the function $\bm{x} \mapsto \tilde{t}(\bm{x})_{s}$ defined in~\eqref{eq:definition-f} is concave on $\left(\R_{+}^{*}\right)^{\X}$.
\end{lemma}
\proof{Proof.}
Fix some $s \in \X$. Concavity is preserved by non-negative weighted sums (section~3.2.1, \cite{boyd2004convex}). Therefore, it is enough to show that the map $\bm{x} \mapsto \min_{\bm{p} \in \U_{sa}}  \prod_{s' \in \X} x_{s'}^{\gamma p_{s'}}$ is concave for each pair $a \in  \A.$ Since concavity is preserved by point-wise minimization (section~3.2.3, \cite{boyd2004convex}), we simply need to show that $\bm{x} \mapsto \prod_{s' \in \X} x_{s'}^{\gamma p_{s'}}$ is concave for each $\bm{p} \in \U_{sa}$. This follows directly from fact 2 in section 3 in \cite{boros2017convex} and the fact that the map $x \mapsto x^{\gamma}$ is monotone and concave since $\gamma \in (0,1)$ and $\bm{x} \ge \bm{0}$. 
\hfill \halmos
\endproof
\begin{example}\label{ex:concavity-tilde-t}
Consider the RMDP instance and the empirical setup of Examples~\ref{ex:non-convex-T} and~\ref{ex:non-convex-tilde-T}. In Figures~\ref{subfig:t-tilde-b-5} and~\ref{subfig:t-tilde-b-10}, we compute $\bm{x}_{1}=\exp_{b}(\bm{v}_{1}),\bm{x}_{2}=\exp_{b}(\bm{v}_{2}),\bm{x}_{\theta} = \theta\bm{x}_{1}+(1-\theta)\bm{x}_{2}$ and  plot the function $\bm{x} \mapsto \tilde{t}_{1}(\bm{x}_{\theta}).$ Note that $\tilde{t}_{1}$ is indeed a concave function (Figures~\ref{subfig:t-tilde-b-5} and~\ref{subfig:t-tilde-b-10}).
\begin{figure}[hbt]
\begin{center}
   \begin{subfigure}{0.4\textwidth}
         \includegraphics[width=1.0\linewidth]{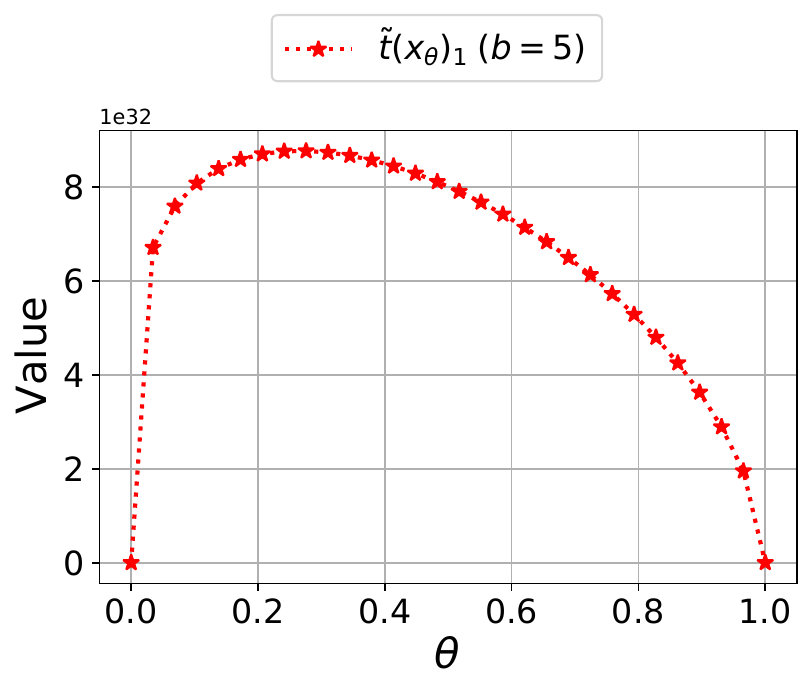}
         \caption{$\tilde{t}$ for $b=5$}
          \label{subfig:t-tilde-b-5}
  \end{subfigure}
     \begin{subfigure}{0.4\textwidth}
         \includegraphics[width=1.0\linewidth]{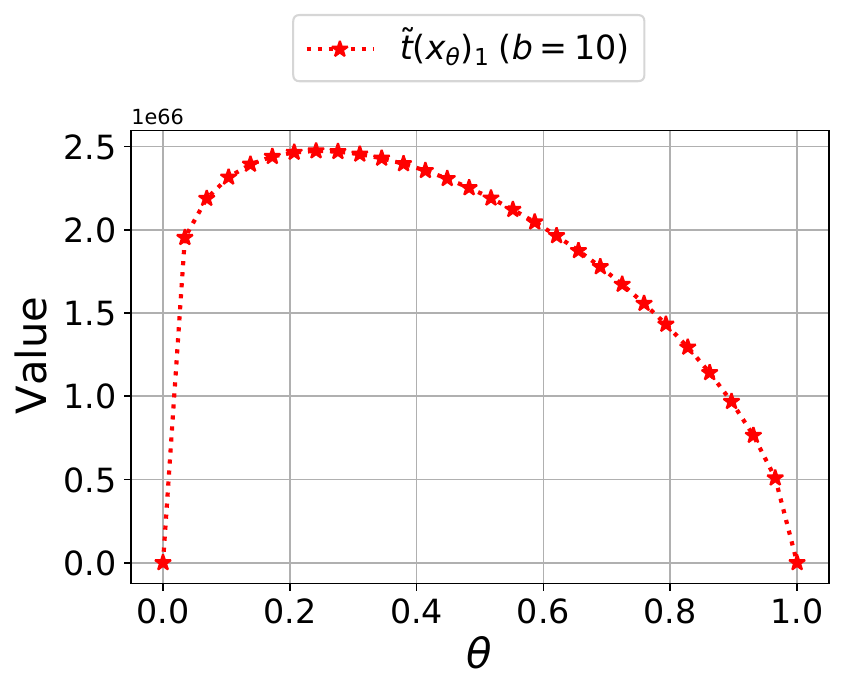}
         \caption{$\tilde{t}$ for $b=10$} \label{subfig:t-tilde-b-10}
  \end{subfigure}
\end{center}
  \caption{Concavity of $\theta \mapsto \tilde{t}_{1}(\bm{x}_{\theta})$ for various values of $b$ (Figures~\ref{subfig:t-tilde-b-5}-\ref{subfig:t-tilde-b-10}).}
  \label{fig:comparing-operators}
\end{figure}
\end{example}
Overall, we obtain the decomposition $\tilde{T} = \log_{b} \circ \tilde{t} \circ \exp_{b}$, with $\tilde{t}$ being a component-wise concave operator. Since $\log_{b} \colon \left(\R_{+}^{*}\right)^{\X} \rightarrow \R^{\X}$ and $\exp_{b}\colon  \R^{\X} \rightarrow \left(\R_{+}^{*}\right)^{\X}$ are one-to-one mappings, we choose the following change of variables:
\[\bm{x} \in \left(\R_{+}^{*}\right)^{\X}, \quad  x_{s}=\exp_{b}(v_{s}), \quad \forall \; s \in \X.\]
We now obtain
\begin{align*}
  \bm{v} \leq \tilde{T}(\bm{v})
  \quad \iff \quad 
  \bm{v} \leq \log_{b} \circ \tilde{t} \circ \exp_{b}(\bm{v}) 
  \quad \iff \quad 
  \exp_{b}(\bm{v}) \leq \tilde{t} \circ \exp_{b}(\bm{v})
  \quad \iff \quad 
  \bm{x} \leq \tilde{t}(\bm{x}),
\end{align*}
and the set $\{ \bm{x} \in \left(\R_{+}^{*}\right)^{\X} \; | \; \bm{x} \leq \tilde{t}(\bm{x})\}$ is a convex set from Lemma~\ref{lem:concavity-of-f}.
\begin{remark}\label{rmk:geometric-programming}
The change of variables $\bm{x} = \exp_{b}(\bm{v}),\bm{v}=\log_{b}(\bm{x})$ is commonly found in the literature on geometric programs~\citep{boyd2007tutorial} and log-log convex programs~\citep{agrawal2019disciplined}. However, the map $\tilde{t}$ is not log-log concave, and the reformulation that we obtain in our main theorem, Theorem~\ref{th:equivalence-programs-v-x}, is not a geometric program, as we highlight in Appendix~\ref{app:geometric-programs}.
\end{remark}
It now remains to choose an objective function $g\colon \R^{\X} \rightarrow \R$ such that $g(\log_{b}(\bm{x}))$ has a  expression. For instance, we can use $g(\bm{v}) = \sum_{s \in \X} \exp_{b}(v_{s})$, which yields $g(\log_{b}(\bm{x}))=\sum_{s \in \X} x_{s}$ after our change of variables. We summarize our results in our main theorem below.
\begin{theorem}\label{th:equivalence-programs-v-x}
Suppose that $\tilde{\bm{v}}\opt$ is the unique fixed point of the regularized robust Bellman operator $\tilde{T}.$ Then $\tilde{\bm{v}}\opt$ is the unique optimal solution to the (potentially) non-convex program
\begin{equation}\label{eq:almost-convex-formulation-regularized-RMDP}
    \max \left\{ \sum_{s \in \X} \exp_{b}(v_{s}) \; | \; \bm{v} \leq \tilde{T}(\bm{v}),\bm{v} \in \R^{\X} \right\}.
\end{equation}
Additionally, $\tilde{\bm{v}}\opt=\log_{b}\left(\tilde{\bm{x}}\opt\right)$, where $\tilde{\bm{x}}\opt \in \R^{\X}$ is the unique optimal solution to the following convex program:
\begin{equation}\label{eq:convex-formulation-regularized-RMDP}
    \max \left\{ \sum_{s \in \X} x_{s} \; | \; 1 \leq x_{s} \leq \tilde{t}(\bm{x})_{s}, \forall \; s \in \X \right\}.
\end{equation}
\end{theorem}
\proof{Proof.}
The optimality and uniqueness of $\tilde{\bm{v}}\opt$ in~\eqref{eq:almost-convex-formulation-regularized-RMDP} follows immediately from the contraction lemma~(Lemma~\ref{lem:contraction-program}). Recall that $\tilde{\bm{v}}\opt \geq \bm{0}$, and therefore we can restrict the domain $\{ \bm{v} \in \R^{\X} \; | \; \bm{v} \leq \tilde{T}(\bm{v})\}$ in~\eqref{eq:almost-convex-formulation-regularized-RMDP} to $\{ \bm{v} \in \R^{\X} \; | \; \bm{v} \leq \tilde{T}(\bm{v}), \bm{v} \geq \bm{0}\}$. We then use our change of variable $\bm{x} = \exp_{b}(\bm{v})$ to obtain that the optimization program~\eqref{eq:almost-convex-formulation-regularized-RMDP} has the same value as~\eqref{eq:convex-formulation-regularized-RMDP}.
The optimization problem in~\eqref{eq:convex-formulation-regularized-RMDP} is convex because its objective function is linear and the feasible set is convex. The feasible set is convex because it is an intersection of two convex sets: the set $\{ \bm{x} \in \R^{\X} \; | x_{s} \geq 1, \forall \; s \in \X\}$, and the set $\{ \bm{x} \in \R^{\X} \; | \; \bm{x} \leq \tilde{t}(\bm{x})\}$. The former set is a half-space, which is well known to be convex, and the latter set is convex by Lemma~\ref{lem:concavity-of-f}. 
\hfill \halmos
\endproof
Solving the convex formulation~\eqref{eq:convex-formulation-regularized-RMDP} to obtain $\tilde{\bm{x}}\opt$ and taking the logarithm gives $\tilde{\bm{v}}\opt$, the fixed-point of $\tilde{T}$, and not $\bm{v}\opt$, the fixed-point of $T$. However, Proposition~\ref{prop:T-T-tilde} bounds the error between with $\tilde{\bm{v}}\opt$ and $\bm{v}\opt$. This argument directly shows that the suboptimality gap of $\tilde{\bm{v}}\opt$ can be chosen to be arbitrarily small, as the following theorem summarizes.
\begin{theorem}\label{th:epsilon-approx-v-v-tilde}
Suppose that $\epsilon >0$ and $b \geq \frac{\log(|\A|)}{\epsilon(1-\gamma)}.$ Let $\tilde{\bm{x}}\opt \in \R^{\X}$ an optimal solution to~\eqref{eq:convex-formulation-regularized-RMDP} and let $\tilde{\bm{v}}\opt = \log_{b}\left( \tilde{\bm{x}}\opt \right)$. Then we have $\tilde{v}\opt_{s}  \leq v\opt_{s} \leq \tilde{v}\opt_{s} + \epsilon$ for each $s \in \X$.
\end{theorem}
\tb{
 We now discuss the {\em main barriers} to the practical use of the optimization program in~\eqref{eq:convex-formulation-regularized-RMDP} for solving RMDPs.
Albeit~\eqref{eq:convex-formulation-regularized-RMDP} is a convex program with $|\X|$ variables and $2 |\X|$ constraints, two main challenges arise in solving it.

The {\em first challenge} in using~\eqref{eq:convex-formulation-regularized-RMDP} is the magnitude of the coefficient $\omega_{b}(s,a) = (1/|\A|) \cdot \exp(b \cdot r_{sa})$ appearing in the expression of $\tilde{t}$ as in~\eqref{eq:definition-f}. These coefficients  may become very large as $b$ increases: in the RMDP instance of Example~\ref{ex:concavity-tilde-t}, $\tilde{t}_{1}$ attains up to $10^{66}$ for $b=10$ (Figure~\ref{fig:comparing-operators}). Using the value $b = \frac{\log(|\A|)}{\epsilon(1-\gamma)}$ from Theorem~\ref{th:epsilon-approx-v-v-tilde}, we obtain that \[\omega_{b}(s,a) = |\A|^{\frac{r_{sa}}{\epsilon (1-\gamma)} -1}\] which is large even for moderate values of $\epsilon>0, \gamma \in (0,1)$ and $r_{sa} \in \R_{+}$. This limitation arises from the change of variables $x_{s} = \exp\left(b \cdot v_{s}\right), \forall \; s \in \X$. At a high-level, if $\bm{v} \in \R^{\X}$ is a value function as defined in~\eqref{eq:value-policy-nominal}, then we can expect that $v_{s}$ is in the order $\frac{1}{1-\gamma}$, so that $x_{s}$ is in the order of $\exp\left(\frac{b}{1-\gamma}\right)$, hence the large magnitude of the constraints coefficients $\omega_{b}(s,a)$.

The {\em second challenge} in using~\eqref{eq:convex-formulation-regularized-RMDP} is the inner minimization over $\bm{p} \in \U_{sa}$ in $\tilde{t}$ which appears in the constraints. This is a common issue in robust optimization, where obtaining tractable robust counterparts of the robust formulations is an important research area~\citep{ben2000robust,yu2015distributionally}. Using strong duality, we will show how to obtain conic reformulations for~\eqref{eq:convex-formulation-regularized-RMDP} (without the inner minimization) in Section~\ref{sec:more-concise-reformulation}.
}
 \tb{ 
\begin{remark}[Some intuition on entropic regularization.]
We comment here on some specific properties of entropic regularization that make it relevant for deriving our convex program.
First, the Bellman operator $T$ is not convex nor concave due to its saddle-point ($\max-\min$) expression. Introducing a KL regularization term for the $\max$ variables is a natural way to replace the maximization program by a closed-form expression. 
The regularized operator $\tilde{T}$ then only involves a single minimization program. However, it is still clear that $\tilde{T}$ may not be convex or concave. As the regularization scalar $b$ approaches $+ \infty$, we know that $\tilde{T}$ approaches $T$, and $T$ is not convex or concave. We resolve this issue with a change of variables. The regularized operator $\tilde{T}$ has a closed-form expression in~\eqref{eq:T-tilde-formulation} thanks to Lemma~\ref{lem:kl-fenchel-conjugate}. Therefore, we can change variables as $\bm{x} = \exp_{b}(\bm{v})$, which yields our convex formulation in~\eqref{eq:convex-formulation-regularized-RMDP}. For other choices of regularizers, e.g.~for the squared $\ell_{2}$-norm, the expression of the regularized Bellman operator is more intricate, and it is not clear what change of variables would yield a convex reformulation~(if any exists).
\end{remark}
}

\tb{
\section{A conic  reformulation based on convex conjugates}\label{sec:more-concise-reformulation}
In this section, we provide a reformulation of the optimization problem ~\eqref{eq:convex-formulation-regularized-RMDP} for the case of convex compact sa-rectangular uncertainty sets described by finitely many convex constraints. In the rest of the paper, we use $m \in \N$ for the number of constraints. In particular, we make the following assumption.
\jgc{
\begin{assumption}\label{assumption-convex}
The set $\U$ is sa-rectangular and there exists $m \in \N$ such that for all $(s,a) \in \X \times \A$, there exist some  proper closed convex functions $g_{sa,1},...,g_{sa,m}\colon \R^{\X} \rightarrow \R \cup \{+\infty\}$ such that
\[
  \U_{sa} = \{ \bm{p} \in \Delta(\X) \; | \; g_{sa,i}(\bm{p})\leq 0, \forall \; i \in \{1,...,m\}\}.
\]
Additionally, Slater's condition holds. That is, for $\cI \subseteq \{1,...,m\}$ the set of indices $i$ such that $g_{sa,i}$ is not affine, there exists a vector $\bm{p}$ in the relative interior of $\Delta(\X)$, such that $g_{sa,i}(\bm{p}) < 0, \forall \; i \in \cI.$
\end{assumption}
For $\bm{p} \in \Delta(\X)$, we will use the notation $\bm{g}_{sa}(\bm{p})$ for the vector $\left(g_{sa,1}(\bm{p}),...,g_{sa,m}(\bm{p})\right) \in \R^{m}$.}
As described in Section~\ref{sec:reformulation-concise}, uncertainty sets satisfying Assumption~\ref{assumption-convex} can model virtually all the practical instances of sa-rectangular RMDPs, e.g. based on $\ell_{\infty},\ell_{1}$ and $\ell_{2}$ norms, or KL divergence.
Under Assumption~\ref{assumption-convex}, our convex formulation~\eqref{eq:convex-formulation-regularized-RMDP} becomes
\begin{equation}\label{eq:reformulation--1}
   \begin{aligned}
    \max \quad  & \sum_{s \in \X} x_{s} \\
    \stc \quad  & x_{s} \leq \sum_{a \in \A} \omega_{b}(s,a)  \min_{\bm{p} \in \Delta(\X),  \bm{g}_{sa}(\bm{p})  \leq \bm{0}}  \prod_{s' \in \X} x_{s'}^{\gamma p_{s'}}, \quad  \forall \; s \in \X,\\
   & x_{s} \geq 1, \forall \; s \in \X.
\end{aligned} 
\end{equation}
We now provide a concise formulation of the constraints in~\eqref{eq:reformulation--1}, based on convex duality. To do so, in Section \ref{sec:tractable counterparts - sa-rec} we convert the minimization problem
\begin{equation}\label{eq:min-program-f-p}
    \min_{\bm{p} \in \Delta(\X), \bm{g}_{sa}(\bm{p})\leq\bm{0}}  \prod_{s' \in \X} x_{s'}^{\gamma p_{s'}}
\end{equation}
for each $s\in \X$ and $a\in \A$ appearing in the constraints in~\eqref{eq:reformulation--1} into a maximization program, using the concept of conjugate and perspective functions. This approach is analogous to the construction of tractable robust counterparts in robust optimization and serves to simplify the constraints~\citep{Ben-Tal2009}. 
 Combining this maximization formulation with the convex formulation~\eqref{eq:reformulation--1}, we obtain a concise convex formulation for RMDPs in Section~\ref{sec:reformulation-concise}. In the same section, we provide a conic formulation for the case of polyhedral uncertainty, and we discuss its numerical complexity. For the sake of conciseness, the conic formulations for ellipsoidal uncertainty and KL-based uncertainty are derived in the appendices (Appendix~\ref{app:sa-rec-other-uncertainty}).
}

\tb{
\subsection{Tractable counterparts of the constraints}\label{sec:tractable counterparts - sa-rec}
In this section, we introduce a maximization formulation for the terms~\eqref{eq:min-program-f-p}, which appear on the right-hand side of the constraints of the optimization problem~\eqref{eq:reformulation--1}. 
\jgc{
We rely on classical tools from the optimization literature, namely conjugate and perspective functions. We define these objects following the lines of \cite{zhen2023unified}. The \emph{convex conjugate function} of $f\colon \R^{\X} \rightarrow \R \cup \{+\infty\}$ is the function $f^*\colon \R^{\X} \rightarrow \R \cup \{+\infty\}$ defined as
\[
  f^*(\bm{y}) = \sup_{\bm{p} \in \R^{\X}} \; \bm{y}^{\top}\bm{p} - f(\bm{p}), \forall \; \bm{y} \in \R^{\X}.
\]
The \emph{perspective function} of $f$ is the function $\R_{+} \times \R^{\X} \rightarrow \R\cup \{+\infty\}$ that maps $(\xi,\bm{x})$ to $\xi f(\bm{x}/\xi)$ if $\xi>0$ and maps $(0,\bm{x})$ to $\sup_{\bm{p} \in \dom(f^{*})} \bm{p}\tr\bm{x}$. We overload notations and write $\xi f(\bm{x}/\xi)$ even when $\xi=0$, with the convention that $0 f(\bm{x}/0) := \sup_{\bm{p} \in \dom(f^{*})} \bm{p}\tr\bm{x}$. It is well-known that the conjugate of a function is always convex, and our definition of perspective functions ensures that the perspective function is proper closed convex whenever the function itself is proper closed convex.
}
We will also use the {\em exponential cone} $\Kexp \subset \R^{3}$, defined as
\[\Kexp := \{ (x_{1},x_{2},x_{3}) \in \R^{3} \; | \; x_{3} \leq x_{2} \log(x_{1}/x_{2}), x_{1},x_{2} >0\} \cup \{ (x_{1},0,x_{3}) \in \R^{3} \; | \; x_{3} \leq 0, x_{1} \geq 0\} .\]
$\Kexp$ is a well-studied convex cone, appearing in various contexts such as robust MDPs based on $\phi$-divergence~\citep{ho2022robust}, Fisher markets~\citep{gao2023infinite}, logistic regression, geometric program and entropy maximization~\citep{chandrasekaran2017relative}. 
With these notations, our main result in this section is the following proposition.
\begin{proposition}\label{prop:tractable-counterpart}
  Let $(s,a) \in \X \times \A$ and $\bm{x} \in [1,+ \infty)^{\X}$. Under Assumption \ref{assumption-convex}, we have,
    \begin{equation}\label{eq:simplifying-min-f-3}
    \begin{aligned}
 \min_{\bm{p} \in \Delta(\X), \bm{g}_{sa}(\bm{p})\leq\bm{0}}  \prod_{s' \in \X} x_{s'}^{\gamma p_{s'}} =   \max \quad  & \min_{s' \in \X} \left\{ z_{s'} + y_{s'} \right\} + \left(\frac{1-\gamma}{\gamma}\right) u + \left(\frac{1+\log(\gamma)}{\gamma}\right)\alpha - \sum_{i=1}^{m} \xi_{i}g_{sa,i}^{*}\left(\frac{\bm{y}^{i}}{\xi_{i}}\right) \\ 
 \stc \quad & (x_{s'},\alpha,z_{s'}) \in \Kexp, \forall \; s' \in \X, \\
 & (1,\alpha,u) \in \Kexp, \\
 & \sum_{i =1}^{m} \bm{y}^{i} = \bm{y}, \\
    & \bm{\xi} \in \R_{+}^{m},\bm{y} \in \R^{\X}, \bm{z} \in \R^{\X}, u\in \R, \alpha \in \R_{+}, \bm{y}^{i} \in \R^{\X}, \forall \; i \in [m].
    \end{aligned}
\end{equation}
\end{proposition}
The rest of this section is devoted to proving Proposition \ref{prop:tractable-counterpart}. We define the map $f\colon \Delta(\X) \rightarrow \R$  as 
\begin{equation}\label{eq:definition-f-p}
    f(\bm{p}) = \prod_{s' \in \X} x_{s'}^{\gamma p_{s'}}, \forall \; \bm{p} \in \Delta(\X).
\end{equation}
Using strong duality in convex programs over sets satisfying Slater's condition (section 5.2.3 in \cite{boyd2004convex}), we have from the definition of the conjugate function that
\begin{align}
  \nonumber \min_{\bm{p} \in \U_{sa}} f(\bm{p}) & = \min_{\bm{p} \in \Delta(\X), \bm{g}_{sa}(\bm{p}) \leq \bm{0}} f(\bm{p})  \\
  \nonumber & = \min_{\bm{p} \in \R^{\X}} f(\bm{p}) + \max_{\bm{\xi} \in \R_{+}^{m},\bm{\mu}\in \R_{+}^{\X},\theta \in \R} \bm{\xi}\tr\bm{g}_{sa}(\bm{p}) - \bm{\mu}\tr\bm{p} + \theta \left(1-\bm{e}\tr\bm{p}\right)  \\
    \nonumber & = \max_{\bm{\xi} \in \R_{+}^{m},\bm{\mu}\in \R_{+}^{\X},\theta \in \R} \theta + \min_{\bm{p} \in \R^{\X}} f(\bm{p}) + \bm{\xi}\tr\bm{g}_{sa}(\bm{p}) - \bm{p}\tr\left(\bm{\mu} + \theta \bm{e}\right)  \\
    \nonumber & = \max_{\bm{\xi} \in \R_{+}^{m},\bm{\mu}\in \R_{+}^{\X},\theta \in \R}   \theta - \max_{\bm{p} \in \R^{\X}} - \left(f + \bm{\xi}\tr\bm{g}_{sa} \right)(\bm{p}) + \bm{p}\tr\left(\bm{\mu} + \theta \bm{e}\right)  \\
   \label{eq:opt-dual-f} & = \max_{\bm{\xi} \in \R_{+}^{m},\bm{\mu}\in \R_{+}^{\X},\theta \in \R} \theta - \left(f + \bm{\xi}\tr\bm{g}_{sa} \right)^* \left(\bm{\mu} + \theta \bm{e}\right)
\end{align}
where we write $\bm{\xi}\tr\bm{g}_{sa}$ for the map $\R^{\X} \rightarrow \R \cup \{+\infty\}$ with $(\bm{\xi}\tr\bm{g}_{sa})(\bm{p}) = \sum_{i=1}^{m} \xi_{i} g_{sa,i}(\bm{p}), \forall \; \bm{p} \in \R^{\X}$.

\jgc{
The following lemma relates the conjugate of a sum of functions to their {\em infimal convolution}. 

\begin{lemma}[Theorem 16.4, \cite{rockafellar1970convex}]\label{lem:infimal convolution}
    Let $f_{1},...,f_{m}:\R^{\X} \rightarrow \R \cup \{ + \infty\}$ be proper closed convex functions such that $\cap_{i=1}^{m} {\sf relint}(\dom(f_{i})) \neq \emptyset$. Then
        \[ \left(\sum_{i=1}^{m} f_{i} \right)^{*}(\bm{y}) = \inf \left\{ \sum_{i=1}^{m} f^{*}_{i}(\bm{y}_{i}) \; | \; \bm{y}_{1},...,\bm{y}_{m} \in \R^{\X},
        \sum_{i=1}^{m} \bm{y}_{i} = \bm{y} \right\}, \forall \; \bm{y} \in \R^{\X}.\]
  \end{lemma}
}
We recall that given a function $g:\R^{\X} \rightarrow \R \cup \{ + \infty\}$ and $\xi \in \R_{+}$, we have $(\xi g)^{*}=\xi g^{*}(\bm{y}/\xi$) (section 3.3.2 in \cite{boyd2004convex}). Combining this with Lemma \ref{lem:infimal convolution} to compute $(f + \bm{\xi}\tr\bm{g}_{sa})^{*}(\bm{y})$, we can simplify the maximization program~\eqref{eq:opt-dual-f} as 
\begin{equation}\label{eq:simplifying-min-f-1}
    \begin{aligned}
    \max \quad  & \theta - f^{*}(\bm{y}_{1}) - \sum_{i=1}^{m} \xi_{i}g_{sa,i}^{*}\left(\frac{\bm{y}^{i}}{\xi_{i}}\right) \\ 
    \stc \quad &  \bm{y}_{1} + \bm{y}_{2} =\bm{\mu} + \theta \bm{e}, \\ 
    & \sum_{i=1}^{m} \bm{y}^{i} = \bm{y}_{2}, \\
    & \bm{\xi} \in \R_{+}^{m},\bm{\mu}\in \R_{+}^{\X},\theta \in \R, \bm{y}_{1} \in \R^{\X}, \bm{y}_{2} \in \R^{\X}, \bm{y}^{i} \in \R^{\X},\forall \; i \in \{1,...,m\}.
    \end{aligned}
\end{equation}
The following proposition, proved in Appendix~\ref{app:proof-tractable-counterpart}, establishes the conjugate function of the function $f\colon \Delta(\X) \to \R$ defined as in~\eqref{eq:definition-f-p}. Here, we use the convention that $\alpha \log(\alpha)=0$ for $\alpha=0$.
\begin{proposition}\label{prop:conjugate-f-tilde}
  Let $\bm{x} \in [1,+\infty)^{\X}$. If $\bm{x}= (1,...,1)$, then
  \[
    f^*(\bm{y}) =
    \begin{cases}
      -1 &\text{if } \bm{y} = \bm{0}, \\
      +\infty & \text{otherwise}.
    \end{cases}            
   \]
When $x_{\bar{s}} >1$ for some $\bar{s} \in \X$, then
\[
 f^*(\bm{y}) =
 \begin{cases}
   \frac{\alpha}{\gamma} \log\left(\frac{\alpha}{\gamma}\right) - \frac{\alpha}{\gamma} &\quad \text{if } \bm{y} = \alpha \log (\bm{x}) \text{ for some } \alpha \in \R_+, \\
   +\infty &\quad \text{otherwise}.
 \end{cases}
\]
\end{proposition}
Substituting in the representation of $f^*$ from Proposition~\ref{prop:conjugate-f-tilde} and removing $\bm{\mu} \in \R^{\X}_{+}$ in \eqref{eq:simplifying-min-f-1}, we obtain  
\begin{equation}\label{eq:simplifying-min-f-2}
\begin{aligned}
    \max \quad  & \theta - \frac{\alpha}{\gamma} \log\left(\frac{\alpha}{\gamma}\right) + \frac{\alpha}{\gamma} - \sum_{i=1}^{m} \xi_{i}g_{sa,i}^{*}\left(\frac{\bm{y}^{i}}{\xi_{i}}\right) \\ 
    \stc \quad &  \alpha \log(\bm{x}) + \bm{y}  \geq \theta \bm{e}, \\ 
    & \sum_{i=1}^{m} \bm{y}^{i} = \bm{y}, \\
    & \bm{\xi} \in \R_{+}^{m},\bm{\mu}\in \R_{+}^{\X},\theta \in \R,  \bm{y} \in \R^{\X}, \bm{y}^{i} \in \R^{\X},\forall \; i \in \{1,...,m\}.
    \end{aligned}
\end{equation}
We can solve for the variable $\theta$, to obtain
\begin{equation}\label{eq:simplifying-min-f-2.1}
    \begin{aligned}
    \max \quad  & \min_{s' \in \X} \{ \alpha\log(x_{s'}) + y_{s'} \} - \frac{\alpha}{\gamma} \log\left(\frac{\alpha}{\gamma}\right) + \frac{\alpha}{\gamma} - \sum_{i=1}^{m} \xi_{i}g_{sa,i}^{*}\left(\frac{\bm{y}^{i}}{\xi_{i}}\right) \\ 
    & \sum_{i=1}^{m} \bm{y}^{i} = \bm{y}, \\
    \stc \quad & \bm{\xi} \in \R_{+}^{m},\bm{y} \in \R^{\X},  \alpha \in \R_{+},\bm{y}^{i} \in \R^{\X},\forall \; i \in \{1,...,m\}.
    \end{aligned}
\end{equation}
It is not straightforward to see that the objective in~\eqref{eq:simplifying-min-f-2.1} is a concave function. We can make this more explicit using the following lemma, which shows that the objective in~\eqref{eq:simplifying-min-f-2.1} can be reformulated using $(\alpha,x) \mapsto \alpha\log(x/\alpha)$, the perspective function of the logarithm, which is a concave function. 
We prove Lemma~\ref{lem:alpha-x-log-analysis} in Appendix~\ref{app:proof-tractable-counterpart}.
\begin{lemma}\label{lem:alpha-x-log-analysis}
Let $\alpha >0, x \geq 1$. Then we have
\begin{align*}
     \alpha\log(x) - \frac{\alpha}{\gamma}\log\left(\frac{\alpha}{\gamma}\right) + \frac{\alpha}{\gamma} = \alpha\log \left(\frac{x}{\alpha}\right) + \left(\frac{1-\gamma}{\gamma}\right)  \alpha\log\left(\frac{1}{\alpha}\right)+ 
 \left(\frac{1+\log(\gamma)}{\gamma}\right)\alpha.
\end{align*}
\end{lemma}
Using Lemma~\ref{lem:alpha-x-log-analysis}, we conclude that \eqref{eq:min-program-f-p} is equal to
\begin{equation}\label{eq:simplifying-min-f-2.2}
    \begin{aligned}
  \max \quad  & \min_{s' \in \X} \left\{ \alpha\log\left(\frac{x_{s'}}{\alpha}\right) + y_{s'} \right\} + \left(\frac{1-\gamma}{\gamma}\right)  \alpha\log\left(\frac{1}{\alpha}\right)+ \left(\frac{1+\log(\gamma)}{\gamma}\right)\alpha - \sum_{i=1}^{m} \xi_{i}g_{sa,i}^{*}\left(\frac{\bm{y}^{i}}{\xi_{i}}\right)  \\ 
  & \sum_{i=1}^{m}\bm{y}^{i} = \bm{y}, \\
    \stc \quad& \bm{\xi} \in \R_{+}^{m},\bm{y} \in \R^{\X},  \alpha \in \R_{+}, \bm{y}^{i} \in \R^{\X},\forall \; i \in \{1,...,m\}.
    \end{aligned}
\end{equation}
\jgc{
To complete the proof of Proposition~\ref{prop:tractable-counterpart}, we use the definition of the exponential cone $\Kexp$. In particular, for $(z,\alpha,x) \in \R \times \R_{+}^{*} \times [1,+\infty)$ and $u \in \R$,
\begin{align*}
     z \leq \alpha\log\left(x/\alpha\right) \quad & \iff \quad \left(x,\alpha,z\right) \in \Kexp \\
     u \leq \alpha\log\left(1/\alpha\right) \quad & \iff \quad \left(1,\alpha,u\right) \in \Kexp
\end{align*}
and these equivalences remain true for $\alpha=0$ using the convention that $\alpha\log(\alpha)=0$ for $\alpha=0$ and by definition of $\Kexp$.
}

\begin{remark}
    When replacing $f^*$ by its expression as in Proposition~\ref{prop:conjugate-f-tilde} to obtain the optimization program~\eqref{eq:simplifying-min-f-2}, we have to distinguish the case where $\bm{x} = (1,...,1)$ and the case where $x_{\bar{s}} >1$ for some $\bar{s} \in \X$. Both these cases lead to the same formulation as in~\eqref{eq:simplifying-min-f-2}, since $\max \{ -\alpha \log(\alpha) + \alpha \; | \; \alpha \geq 0\} = 1$. 
\end{remark}
}

\tb{
\subsection{Concise convex formulation under Assumption~\ref{assumption-convex}}\label{sec:reformulation-concise}
Combining Proposition~\ref{prop:tractable-counterpart} with the formulation~\eqref{eq:reformulation--1}, we arrive at the following theorem.
\begin{theorem}\label{th:cvx-reformulation-general}
  Under Assumption~\ref{assumption-convex}, the set of optimal solutions to the convex program~\eqref{eq:convex-formulation-regularized-RMDP} coincides with the set of optimal solutions to the following convex program:
\begin{equation}\label{eq:reformulation-concise-final}
    \begin{aligned}
    \max \quad  & \sum_{s \in \X} x_{s} \\
    \stc \quad & x_{s} \leq \sum_{a \in \A} \omega_{b}(s,a) \left(y_{sas'}+
    z_{sas'}  + \left(\frac{1-\gamma}{\gamma}\right)  u_{sa} + 
 \left(\frac{1+\log(\gamma)}{\gamma}\right)\alpha_{sa} - \sum_{i =1}^{m} \xi_{sai}g_{sa,i}^{*}\left(\frac{\bm{y}^{i}_{sa}}{\xi_{sai}}\right)\right), \forall \; s,s' \in \X,\\
 & (x_{s'},\alpha_{sa},z_{sas'}) \in \Kexp, \forall \; (s,a,s') \in \X \times \A \times \X, \\
 & (1,\alpha_{sa},u_{sa}) \in \Kexp, \forall \; (s,a) \in \X \times \A, \\
    & \sum_{i =1}^{m} \bm{y}_{sa}^{i} = \bm{y}_{sa}, \forall \; (s,a) \in \X \times \A,\\
   & x_{s} \geq 1, \forall \; s \in \X, \\
    &   \bm{\xi}_{sa}\in \R_{+}^{m},\alpha_{sa} \in \R_{+}, \bm{y}_{sa} \in \R^{\X}, \bm{y}^{i}_{sa} \in \R^{\X}, u_{sa} \in \R, \bm{z}_{sa} \in \R^{\X}, \bm{x} \in \R^{\X}, \forall \; (s,a) \in \X \times \A, \forall \; i \in [m]
\end{aligned}
\end{equation}
\end{theorem}
\proof{Proof of Theorem~\ref{th:cvx-reformulation-general}.}
It remains to check that the constraints  in~\eqref{eq:reformulation-concise-final} are convex, i.e., to prove that the term $-\sum_{i =1}^{m} \xi_{sai}g_{sa,i}^{*}\left(\frac{\bm{y}^{i}_{sa}}{\xi_{sai}}\right)$ is jointly concave in $\left(\bm{y}^{i}_{sa}\right)_{i \in [m]}$ and in $\bm{\xi}_{sa}$. This function is the sum of the negative of the perspective functions of $g^{*}_{sa,i}$, the conjugate function of $g_{sa,i}$, for $i=1,...,m$. The conjugate function is always convex, and the perspective function of a convex function is also convex. The negative of a convex function is concave. Therefore, the constraints  in~\eqref{eq:reformulation-concise-final} are convex constraints. 
\hfill \halmos
\endproof
In the rest of this section, we apply our concise convex formulation~\eqref{eq:reformulation-concise-final} for some of the most common uncertainty sets in the RMDP literature. In particular, we focus on uncertainty sets based on polyhedral functions~\citep{ho2021partial}, ellipsoidal functions~\citep{iyengar2005robust} or Kullback-Leibler divergence~\citep{nilim2005robust}, and we show that in these cases,~\eqref{eq:reformulation-concise-final} is a conic program. In particular, for all these uncertainty sets, the expressions of the convex conjugates $g^{*}_{sa,i}$ and their perspective functions are well-known, as we state in the next proposition, see section 3.3.1 in \cite{boyd2004convex}. \jgc{The expressions of the perspectives functions at $\xi=0$ follow from our convention that $\xi g^{*}(\bm{y}/\xi) := \sup_{\bm{p} \in \dom(g)} \bm{p}\tr\bm{y}$ for $\xi=0$. }
\begin{proposition}\label{prop:example-conjugate}
    \begin{enumerate}
        \item (Affine maps)  Let $g:\R^{\X} \rightarrow \R, \bm{p} \mapsto \bm{a}\tr\bm{p} + b$, for $\bm{a} \in \R^{\X}, b \in \R$. Then $\bm{g}^{*}(\bm{y}) = - b$ if $\bm{y}=\bm{a}$, and $\bm{g}^{*}(\bm{y}) = + \infty$ otherwise. 
        \jgc{
       Let $\bm{y} \in \R^{\X}$. For $\xi>0$, we have $\xi g^{*}(\bm{y}/\xi) = - \xi b$ if $\bm{y} = \xi \bm{a}$, otherwise $\xi g^{*}(\bm{y}/\xi) = +\infty$. For $\xi=0$ we have $\xi g^{*}(\bm{y}/\xi) = 0$ if $\bm{y} = \bm{0}$ and $\xi g^{*}(\bm{y}/\xi) = +\infty$ otherwise.
        }
        \item (Ellipsoidal maps) Let $g:\R^{\X} \rightarrow \R, \bm{p} \mapsto \frac{1}{2} \| \bm{p} \|_{2}^{2}$. Then $g^{*}(\bm{y}) = \frac{1}{2} \| \bm{y} \|_{2}^{2}, \forall \; \bm{y} \in \R^{\X}$.
        \jgc{
       Let $\bm{y} \in \R^{\X}$. For $\xi>0$, we have $\xi g^{*}(\bm{y}/\xi) = (1/\xi)\| \bm{y}\|_{2}^{2}$. For $\xi=0$ we have $\xi g^{*}(\bm{y}/\xi) = 0$ if $\bm{y} = \bm{0}$ and $\xi g^{*}(\bm{y}/\xi) = +\infty$ otherwise.
        }
        \item (Kullback-Leibler divergence). Let $\bm{q} \in \Delta(\X)$ and $g:\Delta(\X) \rightarrow \R \cup \{ + \infty \}, \bm{p} \mapsto \kl(\bm{p},\bm{q})$. Then $g^*(\bm{y}) = \log\left(\sum_{s \in \X} q_{s}\exp(y_{s})\right), \forall \; \bm{y} \in \R^{\X}.$ 
        \jgc{
       Let $\bm{y} \in \R^{\X}$. For $\xi>0$, we have $\xi g^{*}(\bm{y}/\xi) = \xi\log\left(\sum_{s \in \X} q_{s}\exp(y_{s}/\xi)\right)$. For $\xi=0$ we have $\xi g^{*}(\bm{y}/\xi) = \max \{ y_{s} \; | \; s \in \X, q_{s} >0\}$. 
        }
    \end{enumerate}
\end{proposition}
For the sake of conciseness, in the main body of this paper we present in full detail the conic formulation for the case of {\em polyhedral uncertainty}. We provide the reformulations for ellipsoidal uncertainty and KL-based uncertainty in Appendix~\ref{app:sa-rec-other-uncertainty}. 
In the case of polyhedral uncertainty, we have
\begin{equation}\label{eq:uncertainty-polyhedral}
    \U_{sa} = \{ \bm{p} \in \Delta(\X) \; | \; \bm{A}_{sa}\bm{p} \leq \bm{c}_{sa}\}, \forall \; (s,a) \in \X \times \A,
\end{equation}
with $\bm{A}_{sa} \in \R^{m \times \X},\bm{c}_{sa} \in \R^{m}$. Polyhedral uncertainty sets typically describe uncertainty based on $\ell_{\infty}$ norm~\citep{givan1997bounded} and $\ell_{1}$ norm~\citep{ho2018fast}, for which efficient algorithms exist to evaluate the Bellman operator, and which can serve to obtain outer approximation of KL-based uncertainty~\citep{iyengar2005robust}. This type of uncertainty has been used in applications of robust MDPs to healthcare, e.g. \cite{goh2018data} using box uncertainty and \cite{grand2023robustness} using budget of uncertainty. 
Let $\bm{g}_{sa}:\bm{p}\mapsto \bm{A}_{sa}\bm{p} - \bm{c}_{sa}$. We have that $(\bm{\xi}\tr\bm{g}_{sa})^{*}(\bm{y}) = \bm{\xi}\tr\bm{c}_{sa}$ if $\bm{y} = \bm{A}_{sa}\tr\bm{\xi}$, and  $(\bm{\xi}\tr\bm{g}_{sa})^{*}(\bm{y}) = + \infty$ otherwise. Plugging this into~\eqref{eq:reformulation-concise-final}, we obtain the following corollary. 

\begin{corollary}\label{cor:reformulation-concise-polyhedral}
Assume that the uncertainty set $\U$ is sa-rectangular and each set $\U_{sa}$ is polyhedral of the form~\eqref{eq:uncertainty-polyhedral}. The set of optimal solutions to the convex program~\eqref{eq:convex-formulation-regularized-RMDP} coincides with the set of optimal solutions to the following convex program:
\begin{equation}\label{eq:reformulation-concise-polyhedral-conic}
    \begin{aligned}
    \max \quad  & \sum_{s \in \X} x_{s} \\
    \stc \quad & x_{s} \leq \sum_{a \in \A} \omega_{b}(s,a) \left(
    z_{sas'}  + \left(\frac{1-\gamma}{\gamma}\right)  u_{sa} + 
 \left(\frac{1+\log(\gamma)}{\gamma}\right)\alpha_{sa} - \bm{c}_{sa}\tr\bm{\xi}_{sa} + \left(\bm{A}_{sa}\tr\bm{\xi}_{sa}\right)_{s'}\right), \forall \; s,s' \in \X,\\
 & (x_{s'},\alpha_{sa},z_{sas'}) \in \Kexp, \forall \; (s,a,s') \in \X \times \A \times \X, \\
 & (1,\alpha_{sa},u_{sa}) \in \Kexp, \forall \; (s,a) \in \X \times \A, \\
   & x_{s} \geq 1, \forall \; s \in \X, \\
    &   \bm{\xi}_{sa}\in \R_{+}^{m},\alpha_{sa} \in \R_{+},  u_{sa} \in \R, \bm{z}_{sa} \in \R^{\X}, \bm{x} \in \R^{\X}, \forall \; (s,a) \in \X \times \A
\end{aligned}
\end{equation}
\end{corollary}
}
\tb{
\paragraph{Numerical complexity and limitations.} The conic formulation~\eqref{eq:reformulation-concise-polyhedral-conic}
optimizes a linear objective over a feasible set consisting of intersections and Cartesian products of nonnegative orthants and exponential cones. This type of conic program can be solved via primal-dual interior-point methods, with convergence guarantees~\citep{dahl2022primal}. However, in contrast to linear or quadratic programs (e.g. section 4.6 in \cite{ben2001lectures}), to the best of our knowledge, when the domain involves exponential cones, there are no precise bounds on the numerical complexity required to return an $\epsilon$-optimal solution, and deriving such a theoretical bound appears out of the scope of this paper. As discussed at the end Section~\ref{sec:convex-formulation}, the magnitudes of the coefficients $\omega_{b}(s,a):= \frac{1}{|\A|}\exp\left(b \cdot r_{sa}\right)$ appearing in the constraints in~\eqref{eq:reformulation-concise-polyhedral-conic} limit the practical implementation of the iterative algorithms from \cite{dahl2022primal} for solving our conic programs. 
}
\tb{
\section{Convex formulations for s-rectangular RMDPs}\label{sec:s-rectangular}
We now provide our convex formulations for s-rectangular RMDPs. Assume that $\U = \times_{s \in \X } \; \U_{s}$ with $\U_{s} \subseteq \left(\Delta(\X)\right)^{\A}$ is a convex compact set. Since the methods used in this section are similar to the case of sa-rectangular RMDPs, we defer all the proofs for this section to Appendix~\ref{app:proof-s-rectangular}. 

It is known that for s-rectangular RMDPs with convex compact uncertainty set, an optimal policy may be chosen stationary~\citep{wiesemann2013robust}. The robust Bellman operator $T\colon \R^{\X} \rightarrow \R^{\X}$ for s-rectangular RMDPs is defined as
\[
  T(\bm{v})_{s} = \max_{\bm{\pi}_{s} \in \Delta(\A)} \min_{\left(\bm{p}_{a}\right)_{a \in \A} \in \U_{s}} \sum_{a \in \A} \pi_{sa} \left(r_{sa} + \gamma  \bm{p}_{a}^\top \bm{v}\right),  \quad  \forall \; s \in \X.
\]
The regularized Bellman operator becomes, for general choice of regularization policy $\nu \in \Pi$,
\begin{equation}\label{eq:T-tilde-definition-s-rec}
    \tilde{T}(\bm{v})_{s} = \max_{\bm{\pi}_{s} \in \Delta(\A)} \min_{\left(\bm{p}_{a}\right)_{a \in \A}\in \U_{s}} \sum_{a \in \A} \pi_{sa}\left(r_{sa} + \gamma  \bm{p}_{a}^\top \bm{v}\right) - \frac{1}{b} \cdot \kl(\bm{\pi}_{s},\bm{\nu}_{s}), \; \forall \; s \in \X,
  \end{equation}
Similarly to sa-rectangular RMDPs, we choose $\nu_{sa} = 1/A, \forall \; (s,a) \in \X \times \A$ in the rest of this section.
The corresponding operator $\tilde{t}\colon \left(\R_{+}^{*}\right)^{\X} \rightarrow \left(\R_{+}^{*}\right)^{\X}$ becomes, for each state $s \in \X$ and each $\bm{x} \in \left(\R_{+}^{*}\right)^{\X}$,
\begin{equation}\label{eq:t-tilde-s-rec}
    \tilde{t}(\bm{x})_{s}=  \min_{\left(\bm{p}_{a}\right)_{a \in \A}\in \U_{s}} \sum_{a \in \A} \omega_{b}(s,a)\prod_{s' \in \X} x_{s'}^{\gamma p_{as'}}.
\end{equation}
It is straightforward to extend Lemma~\ref{lem:concavity-of-f} to show that this new operator $\tilde{t}$ is component-wise concave, which directly leads to the following theorem, analogous to Theorem~\ref{th:equivalence-programs-v-x} for sa-rectangular uncertainty sets.
\begin{theorem}\label{th:equivalence-programs-v-x-s-rec}
Let $\tilde{\bm{v}}\opt \in \R^{\X}$ be the unique fixed-point of the regularized robust Bellman operator $\tilde{T}$ in~\eqref{eq:T-tilde-definition-s-rec}. 
Then $\tilde{\bm{v}}\opt$ can be computed as $\tilde{\bm{v}}\opt=\log_{b}\left(\tilde{\bm{x}}\opt\right)$, with $\tilde{\bm{x}}\opt \in \R^{\X}$ the unique optimal solution to the following convex program:
\begin{equation}\label{eq:reformulation-s-rec-0}
    \max \left\{ \sum_{s \in \X} x_{s} \; | \; 1 \leq x_{s} \leq \tilde{t}(\bm{x})_{s}, \forall \; s \in \X \right\}
\end{equation}
with $\tilde{t}$ as defined in~\eqref{eq:t-tilde-s-rec}. 
Additionally, let $\epsilon>0$ and $b \geq \log(|\A|)/(\epsilon(1-\gamma))$. Let $\bm{v}\opt$ be the optimal value function of the s-rectangular RMDP. 
Then we have $\tilde{v}\opt_{s}  \leq v\opt_{s} \leq \tilde{v}\opt_{s} + \epsilon$ for each $s \in \X$.
\end{theorem}
To obtain a more explicit optimization program, we remove the minimization in the expression of $\tilde{t}$ using convex conjugate. We use the following assumption.
\jgc{
\begin{assumption}\label{assumption-convex-s-rec}
    The set $\U$ is s-rectangular and there exists $m \in \N$ such that for all $s \in \X $, there exist some proper closed convex functions $g_{s,1},...,g_{s,m}:\R^{\X \times \A} \rightarrow \R$ such that
    \[
      \U_{sa} = \{ \bm{p} \in \Delta(\X)^{\A} | \; g_{s}(\bm{p})\leq 0, \forall \; i \in \{1,...,m\}\}.
    \]
Additionally, Slater's condition holds: for $\cI \subseteq \{1,...,m\}$ the set of indices such that $g_{s,i}$ is not affine, there exists a vector $\bm{p}$ in the relative interior of $\Delta(\X)^{\A}$, such that $g_{s,i}(\bm{p}) <0, \forall \; i \in \cI.$
\end{assumption}
}
Under Assumption~\ref{assumption-convex-s-rec}, we can obtain a maximization program for each component of $\tilde{t}$. In particular, we have the following lemma. We provide the proof in Appendix~\ref{app:proof-s-rectangular}.
\begin{lemma}\label{lem:t-tilde-simplified-s-rec}
    For each $\bm{x} \in [1,+\infty)^{\X}$ and $s \in \X$ the following equality holds:
    \begin{align*}
      \tilde{t}(\bm{x})_{s} = \max \; & - \left(\bm{\xi}_{s}\tr\bm{g}_{s}\right)^{*}(\bm{y}_{s})  +   \sum_{a \in \A} \omega_{b}(s,a) \left( \min_{s' \in \X} \{\alpha_{sa}\log(x_{s'}) + \frac{y_{sas'}}{\omega_{b}(s,a)} \} -\frac{\alpha_{sa}}{\gamma} \log \left(\frac{\alpha_{sa}}{\gamma}\right) + \frac{\alpha_{sa}}{\gamma}  \right) \\
    & \bm{\xi}_{s} \in \R^{m}, \bm{y}_{sa} \in \R^{\X}, \alpha_{sa} \in \R_{+}, \forall \; a \in  \A.
    \end{align*}
\end{lemma}
We can use Lemma~\ref{lem:infimal convolution} to obtain a closed-form expression for $\left(\bm{\xi}_{s}\tr\bm{g}_{s}\right)^{*}(\bm{y}_{s})$. Plugging this back into~\eqref{eq:reformulation-s-rec-0} and using Lemma~\ref{lem:t-tilde-simplified-s-rec}, we obtain the following theorem. 

\begin{theorem}\label{th:cvx-reformulation-general-s-rec}
Under Assumption~\ref{assumption-convex-s-rec},
the set of optimal solutions to the convex program~\eqref{eq:reformulation-s-rec-0} coincides with the set of optimal solutions to the following convex program:
\begin{equation}\label{eq:reformulation-concise-1}
    \begin{aligned}
    \max \quad  & \sum_{s \in \X} x_{s} \\
    \stc \quad & x_{s} \leq - \sum_{i =1}^{m} \xi_{si}g_{s,i}^{*}(\bm{y}^{i}_{s}/\xi_{si}) +  \sum_{a \in 
    \A} y_{sas'}+\omega_{b}(s,a) \left( z_{sas'} + \left(\frac{1-\gamma}{\gamma}\right)u_{sa}+\left(\frac{1+\log(\gamma)}{\gamma}\right)\alpha_{sa} \right), \forall \; s,s' \in \X,\\
    & (x_{s'},\alpha_{sa},z_{sas'}) \in \Kexp, \forall \; (s,a,s') \in \X \times \A \times \X,\\
    & (1,\alpha_{sa},u_{sa}) \in \Kexp,\forall \; (s,a) \in \X \times \A,\\
    & \sum_{i =1}^{m} \bm{y}^{i}_{s} = \bm{y}_{s}, \forall \; s \in \X,\\
    & x_{s} \geq 1, \forall \; s \in \X, \\
    &  \bm{\xi}_{s}\in \R_{+}^{m},\alpha_{sa} \in \R_{+}, \bm{y}_{s} \in \R^{ \A \times \X},\bm{y}_{s}^{i} \in \R^{\A \times \X}, u_{sa} \in \R,\bm{z}_{sa} \in \R^{\X}, \bm{x} \in \R^{\X}, \forall \; (s,a) \in \X \times \A, \forall \; i \in [m].
\end{aligned}
\end{equation}
\end{theorem}
The concavity of~\eqref{eq:reformulation-concise-1} can be proved exactly as the concavity of~\eqref{eq:reformulation-concise-final} for sa-rectangular uncertainty. 
To conclude this section, we now present the application of our concise reformulation~\eqref{eq:reformulation-concise-1} in the case of polyhedral s-rectangular uncertainty sets~\citep{ho2021partial}. For the sake of conciseness, we do not derive reformulations for other common s-rectangular uncertainty sets based on ellipsoidal or KL uncertainty~\citep{grand2021scalable}, noting that the same methods as in Section~\ref{sec:reformulation-concise} apply here. In particular, we make the following assumption.
\begin{assumption}\label{assumption-polyhedral-s-rec} The uncertainty set $\U$ is s-rectangular, and for each $s \in \X$ there exist $\bm{A}_{s} \in \R^{m \times (\A \times \X)},\bm{c}_{s} \in \R^{m}$ for which
    \begin{equation}\label{eq:uncertainty-polyhedral-s-rec}
    \U_{s} = \{ \bm{p} \in \Delta(\X)^{\A} \; | \; \bm{A}_{s}\bm{p} \leq \bm{c}_{s}\}, \forall \; s \in \X.
\end{equation}
\end{assumption}
 Similarly as Corollary~\ref{cor:reformulation-concise-polyhedral} for sa-rectangular RMDPs, we obtain the following reformulation of~\eqref{eq:reformulation-concise-1} under Assumption~\ref{assumption-polyhedral-s-rec}:
 \begin{equation}\label{eq:reformulation-concise-polyhedral-s-rec}
    \begin{aligned}
    \max \quad  & \sum_{s \in \X} x_{s} \\
    \stc \quad & x_{s} \leq - \bm{\xi}_{s}\tr\bm{c}_{s} +  \sum_{a \in 
    \A} \left(\bm{A}_{s}\tr\bm{\xi}_{s}\right)_{as'} +\omega_{b}(s,a) \left( z_{sas'} + \left(\frac{1-\gamma}{\gamma}\right)u_{sa}+\left(\frac{1+\log(\gamma)}{\gamma}\right)\alpha_{sa} \right), \forall \; s,s' \in \X,\\
    & (x_{s'},\alpha_{sa},z_{sas'}) \in \Kexp, \forall \; (s,a,s') \in \X \times \A \times \X,\\
    & (1,\alpha_{sa},u_{sa}) \in \Kexp,\forall \; (s,a) \in \X \times \A,\\
    & \sum_{i =1}^{m} \bm{y}^{i}_{s} = \bm{y}_{s}, \forall \; s \in \X,\\
    & x_{s} \geq 1, \forall \; s \in \X, \\
    &  \bm{\xi}_{s}\in \R_{+}^{m},\alpha_{sa} \in \R_{+}, u_{sa} \in \R,\bm{z}_{sa} \in \R^{\X},  \bm{x} \in \R^{\X}, \forall \; (s,a) \in \X \times \A
\end{aligned}
\end{equation}
}
\section{Conclusion}

We describe the first convex formulation of RMDPs with uncertain transition probabilities, under the assumption that the uncertainty set is rectangular. \tb{The main challenge is to overcome the non-convexity of the robust Bellman operator.} Our results are based on (a) an entropic regularization of the robust Bellman operator and (b) a change of variables introducing the exponential of the components of the value function. \tb{We obtain conic reformulations for classical uncertainty sets based on polyhedral, ellipsoidal or relative entropy-based uncertainty. In these cases, our final conic reformulations involve the exponential cone, the quadratic cone and the non-negative orthant. The main difficulty toward practical implementation of the reformulations obtained in this paper is the exponential size of the coefficients appearing in the constraints of our conic programs, originating from our exponential change of variables.}
Our main contribution is to highlight the problem of finding convex formulations for RMDPs and to provide the first successful attempt (despite the limitation in terms of practical implementation). In particular, our results open novel directions of research of interest for RMDPs. First, it may be possible to obtain exact formulations for RMDPs, not based on regularization, for specific uncertainty sets such as $\ell_{p}$-balls or $\phi$-divergence balls around nominal transitions. Second, other approaches could be based on efficient enumeration of the extreme points of the uncertainty sets, for instance based on $\ell_{1}$ ball uncertainty sets and deterministic transitions probabilities, or more efficient reformulations of the sets $\{\bm{v} \in \R^{\X} \; | \; \bm{v} \leq T(\bm{v})\}$ and $\{\bm{v} \in \R^{\X} \; | \; \bm{v} \geq T(\bm{v})\}$. Finally, it could be of interest to study the case where $\gamma \rightarrow 1$, i.e., the case where the discounted return is replaced with the long-run average return. 
\section*{Acknowledgements}
We would like to thank St{\'e}phane Gaubert for insightful comments on the connections between robust MDPs and stochastic games.
\bibliographystyle{informs2014} 
\bibliography{ref}
\begin{APPENDICES}
\tb{
\section{Details for Section~\ref{sec:preliminaries}}\label{app:detail-preliminaries}
In this section we provide more details on MDPs and RMDPs as introduced in Section~\ref{sec:preliminaries}.
\subsection{Algorithms for MDPs}\label{app:alg-mdps}
\paragraph{Iterative algorithms.}
The value iteration algorithm~\eqref{alg:value-iteration} follows directly from applying Banach's iteration to compute the fixed-point of $T_{\bm{P}}$:
\begin{equation}\label{alg:value-iteration}\tag{VI}
    \bm{v}_{0} \in \R^{\X}, \quad  \bm{v}_{k+1} = T_{\bm{P}}(\bm{v}_{k}), \quad  \forall \; k \in \N.
\end{equation}

The policy iteration algorithm~(PI) construct a sequence $\left(\bm{v}_{k}\right)_{k \in \N}$ as follows: $\bm{v}_{0} \in \R^{\X}$, and for all $k \in \N$,
\begin{align}
\pi_{k} &\text{ such that } T_{\bm{P}}^{\pi_k} (\bm{v}_k) \;=\; T_{\bm{P}}(\bm{v}_k),  \label{eq:PI-policy-improv} \\
\bm{v}_{k+1} & = \bm{v}^{\pi_{k}}.
\label{eq:PI-policy-comp}
\end{align}
In particular, PI alternates between \emph{policy improvement} updates~\eqref{eq:PI-policy-improv}, which compute the policy that attains the $\arg \max$ on each component of the Bellman operator, and \emph{policy evaluation} updates~\eqref{eq:PI-policy-comp}, which computes the value function of the current policy. 

\paragraph{Dual linear program.}
The dual program of~\eqref{eq:linear-program-mdp-primal} provides another linear programming formulation:
\begin{equation}\label{eq:linear-program-mdp-dual}
    \begin{aligned}
        \max_{\bm{\mu}} \quad  &  \sum_{(s,a) \in \X \times \A} \mu_{sa}r_{sa}\\
        \stc \quad &  \sum_{a \in \A} \mu_{sa} = \alpha_{s} + \gamma \sum_{s'\in \X} \sum_{a' \in \A}  P_{s'a's} \mu_{s'a'}, \forall s \in \X, \\
        & \bm{\mu} \in \R^{\X \times \A}_{+}.
    \end{aligned}
  \end{equation}
The dual linear program~\eqref{eq:linear-program-mdp-dual} provides additional structural insights on MDPs. In particular, the decision variables $\bm{\mu} \in \R^{\X \times \A}_{+}$ represent the state-action occupancy frequency:
\[
  \mu_{sa} = \sum_{t=0}^{+\infty} \gamma^{t}\PP \left(s_{t}=s,a_{t}=a\right),
\]
and an optimal policy $\pi\opt$ can be recovered from an optimal solution $\bm{\mu}\opt$ of~\eqref{eq:linear-program-mdp-dual}~\citep{puterman2014markov}. 

}
\subsection{Proof of Lemma~\ref{lem:contraction-program}}\label{app:proof-contraction-lemma}
\proof{Proof.}
We first show that $g(\bm{v}\opt) =  \min \{ g(\bm{v}) \; | \; \bm{v} \geq F(\bm{v}) \}.$ Let $\bm{v} \in \R^{\X}$ such that $\bm{v} \geq F(\bm{v})$. Since $F$ is monotone, this shows that $\bm{v} \geq F^{k}(\bm{v}), \forall \; k \geq 0$. But we have $F^{k}(\bm{v}) \rightarrow \bm{v}\opt$, from which we conclude that
for any $\bm{v} \in \R^{\X},$ we have  $\bm{v} \geq F(\bm{v}) \Rightarrow \bm{v} \geq \bm{v}\opt$. Since $g\colon \R^{\X} \rightarrow \R$ is component-wise non-decreasing, this shows that $g(\bm{v}) \geq g(\bm{v}\opt)$. Because $\bm{v}\opt \geq F(\bm{v}\opt) = \bm{v}\opt$, we can conclude that $g(\bm{v}\opt) =  \min \{ g(\bm{v}) \; | \; \bm{v} \geq F(\bm{v}) \}.$
We now show that if $g$ is a component-wise increasing function, $\bm{v}^*$ is the unique vector attaining the $\arg \min$ in $\min \{ g(\bm{v}) \; | \; \bm{v} \geq F(\bm{v}) \}$. Assume that $\bm{v}' \in \R^{\X}$ is such that $\bm{v}'\geq F(\bm{v}')$ and $g(\bm{v}') = g(\bm{v}\opt)$. From $\bm{v}'\geq F(\bm{v}')$ we know that $\bm{v}' \geq \bm{v}\opt$. But if there exists $v_{s}' > v\opt_{s}$ for some $s \in \X$, then $g(\bm{v}') > g(\bm{v}\opt)$ since $g$ is component-wise increasing. This shows that $\bm{v}' = \bm{v}\opt$.

The proof for $g(\bm{v}\opt) =  \max \{ g(\bm{v}) \; | \; \bm{v} \leq F(\bm{v}) \}$ is similar and we omit it for conciseness.
\hfill \halmos
\endproof

\section{Proof of Proposition~\ref{prop:T-T-tilde}}\label{app:prf:prop:T-T-tilde}
\proof{Proof.}
\begin{enumerate}
    \item $\tilde{T}$ is monotone since $\U_{sa} \subseteq \Delta(\X) \subseteq \R_{+}^{\X}$, and $\pi_{sa} \geq 0,$ for all $(s,a) \in \X \times \A$ for $\pi \in \Pi$. \tb{Since we assume that $\nu$ is the uniform policy, we note that $\kl(\bm{\pi}_{s},\bm{\nu}_{s}) < +\infty$ for any $\bm{\pi}_{s} \in \Delta(\A)$.} Let $\bm{v}_{1},\bm{v}_{2} \in \R^{\X}$ and $s \in \X$. Then 
    \tb{
    \begin{align*}
        \tilde{T}(\bm{v}_{1})_{s} & = \max_{\bm{\pi}_{s} \in \Delta(\A)}  \sum_{a \in \A} \pi_{sa}\left(r_{sa} + \gamma \min_{\bm{p} \in \U_{sa}} \bm{p}^\top \bm{v}_{1}\right) - (1/b)\cdot \kl(\bm{\pi}_{s},\bm{\nu}_{s}) \\
        & = \max_{\bm{\pi}_{s} \in \Delta(\A)}  \sum_{a \in \A} \pi_{sa}\left(r_{sa} + \gamma \min_{\bm{p} \in \U_{sa}} \bm{p}^\top ( \bm{v}_{2} + \bm{v}_{1}-\bm{v}_{2})\right) - (1/b)\cdot \kl(\bm{\pi}_{s},\bm{\nu}_{s})\\
         & \leq \max_{\bm{\pi}_{s} \in \Delta(\A)}  \sum_{a \in \A} \pi_{sa}\left(r_{sa} + \gamma \min_{\bm{p} \in \U_{sa}} \bm{p}^\top ( \bm{v}_{2} + \|\bm{v}_{1}-\bm{v}_{2}\|_{\infty} \bm{e})\right) - (1/b)\cdot \kl(\bm{\pi}_{s},\bm{\nu}_{s})  \\
         & = \max_{\bm{\pi}_{s} \in \Delta(\A)}  \sum_{a \in \A} \pi_{sa}\left(r_{sa} + \gamma \min_{\bm{p} \in \U_{sa}} \bm{p}^\top \bm{v}_{2} \right) - (1/b)\cdot \kl(\bm{\pi}_{s},\bm{\nu}_{s}) + \gamma \|\bm{v}_{1}-\bm{v}_{2}\|_{\infty}  \\
        & = \tilde{T}(\bm{v}_{2})_{s} + \gamma \| \bm{v}_{1}-\bm{v}_{2} \|_{\infty},
    \end{align*}
    }
    since $\sum_{s' \in \X} p_{s'}=1,\forall \; \bm{p} \in \U_{sa}$, and $\sum_{a \in \A} \pi_{sa}=1, \forall \; \bm{\pi}_{s} \in \Delta(\A)$.
    Similarly, we can show that $\tilde{T}(\bm{v}_{2})_{s}\leq \tilde{T}(\bm{v}_{1})_{s} + \gamma \| \bm{v}_{1}-\bm{v}_{2} \|_{\infty},$
    and therefore, we conclude that $\tilde{T}$ is a contraction for the $\ell_{\infty}$ norm:
    \[ \| \tilde{T}(\bm{v}_{1})-\tilde{T}(\bm{v}_{2})\|_{\infty} \leq \gamma \| \bm{v}_{1}-\bm{v}_{2} \|_{\infty}.\]
    \item \tb{We now use the fact that $\nu$ is the uniform policy. In this case,} we have $0 \leq (1/b)\cdot \kl(\bm{\pi}_{s},\bm{\nu}_{s}) \leq \log(|\A|)/b$ for any $\bm{\pi}_{s} \in \Delta(\A)$. Therefore, we have
    \begin{align*}
        \tilde{T}(\bm{v})_{s} & = \max_{\bm{\pi}_{s} \in \Delta(\A)}  \sum_{a \in \A} \pi_{sa}\left(r_{sa} + \gamma \min_{\bm{p} \in \U_{sa}} \bm{p}^\top \bm{v}\right) - (1/b)\cdot \kl(\bm{\pi}_{s},\bm{\nu}_{s}) \\& \geq \max_{\bm{\pi}_{s} \in \Delta(\A)}  \sum_{a \in \A} \pi_{sa}\left(r_{sa} + \gamma \min_{\bm{p} \in \U_{sa}} \bm{p}^\top \bm{v}\right) - \log(|\A|)/b, \\
        \tilde{T}(\bm{v})_{s} & = \max_{\bm{\pi}_{s} \in \Delta(\A)}  \sum_{a \in \A} \pi_{sa}\left(r_{sa} + \gamma \min_{\bm{p} \in \U_{sa}} \bm{p}^\top \bm{v}\right) - (1/b)\cdot \kl(\bm{\pi}_{s},\bm{\nu}_{s}) \\
        &\leq \max_{\bm{\pi}_{s} \in \Delta(\A)}  \sum_{a \in \A} \pi_{sa}\left(r_{sa} + \gamma \min_{\bm{p} \in \U_{sa}} \bm{p}^\top \bm{v}\right).
    \end{align*}
    It is straightforward to conclude since
    \[T(\bm{v})_{s} = \max_{\bm{\pi}_{s} \in \Delta(\A)}  \sum_{a \in \A} \pi_{sa}\left(r_{sa} + \gamma \min_{\bm{p} \in \U_{sa}} \bm{p}^\top \bm{v}\right).\]
    \item Let $\bm{v} \in \R^{\X}$. We have shown that $\tilde{T}(\bm{v})_{s} \leq T(\bm{v})_{s} \leq \tilde{T}(\bm{v}) + \log(|\A|)/b$, for $s \in \X$. Let us apply this with $\bm{v} = \tilde{\bm{v}}\opt$, the unique fixed-point of $\tilde{T}$. We obtain
    \[  \tilde{v}\opt_{s} \leq T(\tilde{\bm{v}}\opt)_{s} \leq \tilde{v}\opt_{s}+ \log(|\A|)/b.\]
    If we repeatedly apply the operator $T$ to both sides of the inequality $\tilde{\bm{v}}\opt \leq T(\tilde{\bm{v}}\opt)$, we obtain
    \[\tilde{\bm{v}}\opt \leq T^{\ell}(\tilde{\bm{v}}\opt), \forall \; \ell \in \N.\]
    But $T^{\ell}(\tilde{\bm{v}}\opt) \rightarrow \bm{v}\opt$, where $\bm{v}\opt$ is the unique fixed-point of $T$. Therefore, $\tilde{\bm{v}}\opt \leq \bm{v}\opt$.
     Now if we repeatedly apply the operator $T$ to both sides of the inequalities $T(\tilde{\bm{v}}\opt)_{s} \leq \tilde{v}\opt_{s} + \log(|\A|)/b$ for all $s \in \X$, we obtain 
     \[T^{\ell}(\tilde{\bm{v}}\opt)_{s} \leq \tilde{v}\opt_{s} + \frac{\log(|\A|)}{b}\left( \sum_{k=0}^{\ell-1} \gamma^{k}\right).\]
     Taking the limit when $\ell \rightarrow + \infty$, we obtain that $v\opt_{s} \leq \tilde{v}\opt_{s} + \frac{\log(|\A|)}{b \cdot (1-\gamma)}$.
\end{enumerate}
\hfill \halmos
\endproof

\section{Relation with geometric and log-log convex programs}\label{app:geometric-programs}
\paragraph{Geometric programs.}
Geometric programs~\citep{boyd2007tutorial} are optimization programs of the form
\begin{equation}\label{eq:geometric-program}
\begin{aligned}
    \min & \; h_{0}(\bm{x}) \\
    & h_{i}(\bm{x}) \leq 1, \forall \; i \in \{1,...,m\},\\
    & k_{i}(\bm{x}) =1, \forall \; i \in \{1,...,p\}, \\
    & \bm{x} \in \left(\R_{+}^{*}\right)^{n},
\end{aligned}
\end{equation}
where $k_{i}\colon \left(\R_{+}^{*}\right)^{n} \rightarrow \R$ are monomials:
\[k_{i}(\bm{x}) = c_{i} \prod_{j=1}^{n} x_{j}^{a_{ij}}, c_{i}>0,\left(a_{ij}\right)_{j} \in \R^{\X},\]
and $h_{i}\colon \left(\R_{+}^{*}\right)^{n} \rightarrow \R$ are sums of monomials. Under the change of variables $v_{j}=\log(x_{j})$, geometric programs with variable $\bm{x} \in \left(\R_{+}^{*}\right)^{n}$ can be reformulated as equivalent convex programs with a variable $\bm{v} \in \R^{\X}$. In particular,
\[ k_{i}(\bm{x}) = k_{i}(\exp(\bm{v})) = \exp\left(\log(c_{i}) + \sum_{j=1}^{n} a_{ij}v_{j}\right),\]
and the constraints $k_{i}(\bm{x}) \leq 1$, reformulated as $\log(k_{i}(\bm{x}))\leq 0$, becomes linear in the variable $\bm{v}$, while the constraints $h_{i}(\bm{x}) \leq 1$ can be reformulated using the log-sum-exp function, which is convex. Our formulation~\eqref{eq:convex-formulation-regularized-RMDP} resembles the geometric program~\eqref{eq:geometric-program}, except for two crucial differences: the direction of the inequalities in the constraints are reversed, and the presence of the minimization term over $\U_{sa}$ in the expression of $\bm{x} \mapsto \tilde{t}(\bm{x})_{s}$ as defined in~\eqref{eq:definition-f}.
In particular, if we attempt to reformulate~\eqref{eq:convex-formulation-regularized-RMDP} as a geometric program, we see that ~\eqref{eq:convex-formulation-regularized-RMDP} is equal to
\begin{equation}\label{eq:geo-formulation-regularized-RMDP}
    \begin{aligned}
            \min & \; - \sum_{s \in \X} x_{s} \\
            x_{s} + & \max_{\left(\bm{p}_{a}\right)_{a \in \A} \in \times_{a \in \A} \U_{sa} } - \sum_{a \in \A} \omega_{b}(s,a)   \prod_{s' \in \X} x_{s'}^{\gamma p_{as'}} \leq 0, \; \forall \; s \in \X,\\
            x_{s} & \geq 1, \forall \; s \in \X.
    \end{aligned}
\end{equation}
Because of the $\max$ term in each constraint,~\eqref{eq:geo-formulation-regularized-RMDP} is not a geometric program, and because of the negative term $-\frac{1}{|\A|}\exp_{b}(r_{sa})$ in each of the terms $-\frac{1}{|\A|}\exp_{b}(r_{sa})\prod_{s' \in \X} x_{s'}^{\gamma p_{s'}}$, the optimization program~\eqref{eq:geo-formulation-regularized-RMDP} is not a {\em robust} geometric program~\citep{hsiung2008tractable}.

\paragraph{Log-log convex functions.}
Another important class of problems that resemble~\eqref{eq:convex-formulation-regularized-RMDP} is the class of {\em log-log convex programs}~\citep{agrawal2019disciplined}, which generalizes geometric programs. Log-log convex programs maximize {\em log-log concave functions} under log-log convex constraints. In particular, a function $f\colon \R^{\X} \rightarrow \R_{+}^{*}$ is log-log concave if its logarithm is concave in the exponent of its arguments, i.e.,  if $F=\log \circ f \circ \exp$ is concave. Our functions $\bm{x} \mapsto \tilde{t}(\bm{x})_{s}$ as defined in~\eqref{eq:definition-f} are not log-log concave, since because $\log \circ \tilde{t} \circ \exp$ is equal to the operator $\tilde{T}$ as defined in~\eqref{eq:T-tilde-formulation}, which may be neither component-wise convex nor concave, as highlighted in Example~\ref{ex:non-convex-tilde-T}.

\tb{
\section{Proof for Section \ref{sec:tractable counterparts - sa-rec}}\label{app:proof-tractable-counterpart}
\proof{Proof of Proposition \ref{prop:conjugate-f-tilde}}
We first treat the case where $\bm{x}=(1,...,1).$ In this case, algebraic manipulation shows that $f(\bm{p}) = 1, \forall \; \bm{p} \in \R^{\X}$ with $f$ defined in~\eqref{eq:definition-f-p}. By the definition of the convex conjugate $f^*$, we have that $f^*(\bm{y})= \sup_{\bm{p} \in \R^{\X}} \bm{y}\tr\bm{p} -1$. Therefore, $f^*(\bm{y})=+\infty$ if $\bm{y} \neq \bm{0}$, and $f^*(\bm{y})=-1$ if $\bm{y} = \bm{0}$.

We now treat the case where $\exists \; \bar{s} \in \X, x_{\bar{s}} >1$. Since $\bm{p} \mapsto \log \circ f(\bm{p})$ is an affine function, it is easy to compute the level sets of $f$, i.e., the sets $\cL_{\alpha} = \{ \bm{p} \in \R^{\X} \; | \; f(\bm{p}) = \alpha\}$ for $\alpha >0$. This will prove helpful in computing $f^*$. In particular, we have
\[\bm{p} \in \cL_{\alpha} \iff f(\bm{p}) = \alpha\iff \log(f(\bm{p})) = \log(\alpha)  \iff \bm{p}\tr \log(\bm{x}) = \log(\alpha)/\gamma.\]
This shows that $\cL_{\alpha}$ is a non-empty hyperplane with normal vector $\log(\bm{x})$, with $\log(\bm{x}) \neq 
\bm{0}$, and $\cL_{\alpha}$ can be written as
\[\cL_{\alpha} = \{ \bm{p} \in \R^{\X} \; | \; \bm{p}\tr \log(\bm{x}) = \log(\alpha)/\gamma\}. \]
We now have, for $\bm{y} \in \R^{\X}$,
\[    f^*(\bm{y}) = \sup_{\bm{p} \in \R^{\X}} \; \bm{y}^{\top}\bm{p} - f(\bm{p})  \geq \sup_{\bm{p} \in \cL_{\alpha}} \; \bm{y}^{\top}\bm{p} - f(\bm{p})  = \sup_{\bm{p} \in \cL_{\alpha}} \; \bm{y}^{\top}\bm{p} - \alpha.\]
Since $\cL_{\alpha}$ is a hyperplane with normal vector $\log(\bm{x})$,  either $\bm{y} = \beta \log(\bm{x})$ for some $\beta \in \R$ or $\sup_{\bm{p} \in \cL_{\alpha}} \; \bm{y}^{\top}\bm{p} = + \infty.$
Therefore, we have shown that  
$f^*(\bm{y}) < + \infty \Rightarrow \bm{y} \in \{ \beta \log(\bm{x}) \; | \; \beta \in \R\}$. Let $\beta \in \R$ and  $\bm{y}=\beta \log(\bm{x})$, we have also shown that
\[
  f^*(\bm{y}) \quad\geq \quad\beta \frac{\log(\alpha)}{\gamma} - \alpha, \quad \forall \; \alpha >0.
\]
Therefore, if $\beta<0$, taking the limit as $\alpha \rightarrow 0$ in the above inequality we have $f^*(\bm{y}) = + \infty$ and we have shown that $f^*(\bm{y}) < + \infty \Rightarrow \bm{y} \in \{ \alpha \log(\bm{x}) \; | \; \alpha \geq 0\}.$

Next, we show that $  \bm{y} \in \{ \alpha \log(\bm{x}) \; | \; \alpha \geq 0\} \Rightarrow f^*(\bm{y}) < + \infty.$ Let $\bm{y} = \alpha \log(\bm{x})$ with $\alpha \geq 0.$ If $\alpha=0$, then $f^*(\bm{y})=f^*(\bm{0}) = \sup_{\bm{p} \in \R^\X} - f(\bm{p}) = 0$. We now treat the case $\alpha>0$. Setting the gradient of $\bm{p} \mapsto \bm{y}^{\top}\bm{p} - f(\bm{p})$ to $\bm{0}$ gives
$y_{s} - \gamma \log(x_{s})f(\bm{p}) = 0, \forall \; s \in \X.$ Since $\bm{y}=\alpha \log(\bm{x})$, we conclude that
\[
  (\alpha - \gamma f(\bm{p}))\log(x_{s}) = 0, \forall \; s \in \X.
\] 
Since $x_{\bar{s}} \neq 1$, we have $\log(x_{\bar{s}}) \neq 0$, and therefore $\alpha = \gamma f(\bm{p})$.
We can conclude that the maximum of $\bm{p} \mapsto \bm{y}^{\top}\bm{p} - f(\bm{p})$ is attained at (any) $\bm{p} \in \cL_{\alpha/\gamma}$ and $f^*(\bm{y}) = (\alpha \log(\alpha/\gamma) - \alpha)/\gamma.$ Overall, we have proved that
    \[
 f^*(\bm{y}) =
 \begin{cases}
   \frac{\alpha}{\gamma} \log\left(\frac{\alpha}{\gamma}\right) - \frac{\alpha}{\gamma} &\quad \text{if } \bm{y} = \alpha \log (\bm{x}) \text{ for some } \alpha \in \R_+, \\
   +\infty &\quad \text{otherwise}.
 \end{cases}
\]
\hfill \Halmos \endproof

\proof{Proof of Lemma~\ref{lem:alpha-x-log-analysis}.}
Let $\alpha>0,x\geq 1$. We have
\begin{align*}
    \alpha\log(x) - \frac{\alpha}{\gamma}\log\left(\frac{\alpha}{\gamma}\right) + \frac{\alpha}{\gamma} & = \alpha\log(x) - \alpha \log(\alpha) + \alpha\log(\alpha) - \frac{\alpha}{\gamma}\log\left(\frac{\alpha}{\gamma}\right) + \frac{\alpha}{\gamma} \\
    & = \alpha\log \left(\frac{x}{\alpha}\right) + \alpha\log(\alpha) - \frac{\alpha}{\gamma}\log\left(\frac{\alpha}{\gamma}\right) + \frac{\alpha}{\gamma}.
\end{align*}
We now turn to reformulating $\alpha\log(\alpha) - \frac{\alpha}{\gamma}\log\left(\frac{\alpha}{\gamma}\right) + \frac{\alpha}{\gamma}.$ We have
\begin{align*}
    \alpha\log(\alpha) - \frac{\alpha}{\gamma}\log\left(\frac{\alpha}{\gamma}\right) + \frac{\alpha}{\gamma} & = \alpha\log(\alpha) - \frac{\alpha}{\gamma}\log\left(\alpha\right) + \frac{\alpha}{\gamma}\log\left(\gamma\right) + \frac{\alpha}{\gamma} \\
    & = \alpha\log(\alpha)\left(\frac{\gamma-1}{\gamma}\right) + \frac{\alpha}{\gamma}\log(\gamma) + \frac{\alpha}{\gamma}\\
    & = \left(\frac{1-\gamma}{\gamma}\right)  \alpha\log\left(\frac{1}{\alpha}\right)+ \frac{\alpha}{\gamma}\left(1+\log(\gamma)\right).
\end{align*}
Overall, we have shown that
\[  \alpha\log(x) - \frac{\alpha}{\gamma}\log\left(\frac{\alpha}{\gamma}\right) + \frac{\alpha}{\gamma} = \alpha\log \left(\frac{x}{\alpha}\right) + \left(\frac{1-\gamma}{\gamma}\right)  \alpha\log\left(\frac{1}{\alpha}\right)+ \frac{\alpha}{\gamma}\left(1+\log(\gamma)\right).\]
\hfill \halmos
\endproof
}

\tb{

\section{Conic reformulation for some classical sa-rectangular uncertainty sets}\label{app:sa-rec-other-uncertainty}
In this appendix, we derive the conic formulations for the case of ellipsoidal uncertainty and KL-based uncertainty. 
\subsection{Ellispoidal uncertainty}
We start by introducing ellipsoidal uncertainty.
Assume that
\begin{equation}\label{eq:uncertainty-ellipsoidal}
    \U_{sa} = \left\{ \bm{p} \in \Delta(\X) \; | \; \frac{1}{2}\| \bm{p} - \bm{P}^{0}_{sa} \|_{2}^{2} \leq \theta_{sa}\right\}, \forall \; (s,a) \in \X \times \A,
\end{equation}
with $\bm{P}^{0}$ the nominal transition probabilities and $\theta_{sa} \geq 0$. This type of uncertainty sets naturally arise when using a second-order approximation of the log-likehood function~\citep{nilim2005robust}. Let $g_{sa}:\bm{p} \mapsto \frac{1}{2}\| \bm{p} - \bm{P}^{0}_{sa} \|_{2}^{2} - \theta_{sa}$. Proposition~\ref{prop:example-conjugate} gives that $g_{sa}^{*}(\bm{y}) = \theta_{sa}  + \bm{y}\tr\bm{P}^{0}_{sa} + \frac{1}{2} \| \bm{y}\|_{2}^{2}$ for any $\bm{y} \in \R^{\X}$. 
\jgc{
For $\xi > 0$ we have $\xi g_{sa}^{*}(\bm{y}/\xi) = \xi \theta_{sa} + \bm{y}\tr\bm{P}^{0}_{sa} + \frac{1}{2 \xi} \| \bm{y}\|_{2}^{2}$, and for $\xi=0$ we have that $\xi g_{sa}^{*}(\bm{y}/\xi) = 0$ if $\bm{y}=\bm{0}$ and $\xi g_{sa}^{*}(\bm{y}/\xi) = +\infty$ otherwise.
}
To obtain a conic program, we will make use of the {\bf $n$-dimensional rotated quadratic cone} $\cQ^{n}_{r}$, defined as \[\cQ^{n}_{r} := \{ \bm{y} \in \R^{n} \; | \; \sum_{s=3}^{n} y_{s}^{2} \leq 2 y_{1} y_{2}\}.\] 
Indeed, following the definition of $\cQ^{n}_{r}$, we directly obtain, for $(\xi,q,\bm{y}) \in \R_{+}^{*} \times \R \times \R^{n}$,
\[ \frac{1}{2\xi} \| \bm{y}\|_{2}^{2} \leq q \iff (\xi,q,\bm{y}) \in \cQ^{n+2}_{r}.\]
\jgc{Note that this equivalence remains true at $\xi=0$.}
The rotated quadratic cone $\cQ^{n}_{r}$ can be reformulated as the quadratic cone $\cQ^{n} := \{ \bm{y} \in \R^{n} \; | \; \sum_{s=2}^{n} y_{s}^{2} \leq  y_{1}^{2}\}$ by introducing auxiliary variables~\citep{mosek}. Optimizing over quadratic cones, also known as {\em second-order cone optimization}, is one of the most studied areas of conic optimization, see lecture 4 in \cite{ben2001lectures}. 
Plugging this into~\eqref{eq:reformulation-concise-final}, we obtain the following corollary. 
\begin{corollary}\label{cor:reformulation-concise-ellipsoidal}
Assume that the uncertainty set $\U$ is sa-rectangular and each set $\U_{sa}$ is ellipsoidal of the form~\eqref{eq:uncertainty-ellipsoidal}. Then the convex formulation~\eqref{eq:reformulation-concise-final} can be reformulated as follows:
\begin{equation}\label{eq:reformulation-concise-ellipsoid-conic}
    \begin{aligned}
    \max \quad  & \sum_{s \in \X} x_{s} \\
    \stc \quad & x_{s} \leq \sum_{a \in \A} \omega_{b}(s,a) \left(
    y_{sas'} + z_{sas'}  + \left(\frac{1-\gamma}{\gamma}\right)  u_{sa} + 
 \left(\frac{1+\log(\gamma)}{\gamma}\right)\alpha_{sa} - \xi_{sa}\theta_{sa} - \bm{y}_{sa}\tr\bm{P}^{0}_{sa} - q_{sa}\right), \forall \; s,s' \in \X,\\
 & (x_{s'},\alpha_{sa},z_{sas'}) \in \Kexp, \forall \; (s,a,s') \in \X \times \A \times \X, \\
 & (1,\alpha_{sa},u_{sa}) \in \Kexp, \forall \; (s,a) \in \X \times \A, \\
 & (\xi_{sa},q_{sa},\bm{y}_{sa}) \in \cQ^{|\X|+2}_{r}, \forall \; (s,a) \in \X \times \A, \\
   & x_{s} \geq 1, \forall \; s \in \X, \\
    &   \bm{\xi}_{sa}\in \R_{+}^{m},\alpha_{sa} \in \R_{+},  u_{sa} \in \R, \bm{z}_{sa} \in \R^{\X},\bm{y}_{sa} \in \R^{\X}, q_{sa} \in \R,  \bm{x} \in \R^{\X}, \forall \; (s,a) \in \X \times \A.
\end{aligned}
\end{equation}
\end{corollary}
\subsection{Uncertainty based on KL divergence}
Models based on Kullback-Leibler divergence can be written \begin{equation}\label{eq:uncertainty-kl}
    \U_{sa} = \{ \bm{p} \in \Delta(\X) \; | \; \kl\left(\bm{p},\bm{P}^{0}_{sa}\right) \leq \theta_{sa}\}, \forall \; (s,a) \in \X \times \A,
\end{equation}
with $\bm{P}^{0}$ the nominal transition probabilities and $\theta_{sa} \geq 0$. This model is popular both due its statistical guarantees and its tractability, see section 5 in \cite{iyengar2005robust}. We define $g_{sa}:\bm{p} \mapsto \kl\left(\bm{p},\bm{P}^{0}_{sa}\right) - \theta_{sa}$. Proposition~\ref{prop:example-conjugate} gives that $g_{sa}^{*}(\bm{y}) = \theta_{sa} + \log\left(\sum_{s' \in \X} P^{0}_{sas'}\exp(y_{s'}) \right)$ for $\bm{y} \in \R^{\X}$, so that for $\xi > 0$ we have $\xi g_{sa})^{*}(\bm{y}/\xi) = \xi\theta_{sa} + \xi\log\left(\sum_{s' \in \X} P^{0}_{sas'}\exp(y_{s'}/\xi) \right)$ for any $\bm{y} \in \R^{\X}$. 
\jgc{
For $\xi=0$, our definition of perspectives functions at $\xi=0$ ensures that $\xi g_{sa}^{*}(\bm{y}/\xi) = \max \{ y_{s'} \; | \; s' \in \X, P_{sas'}^{0} >0\}$.
}
From chapter 7.1.1 in \cite{mosek}, we know that, for $(\xi,q,\bm{y}) \in \R_{+}^{*} \times \R \times \R^{n}$, we have
\[\xi\log\left(\sum_{s' \in \X} P^{0}_{sas'}\exp(y_{s'}/\xi)\right) \leq q \iff
\begin{cases}
   & \sum_{s' \in \X} \kappa_{s'} \leq \xi, \\
    & (\kappa_{s'}/P_{sas'}^{0}, \xi,y_{s'} - q)  \in \Kexp, \forall \; s' \in \X \; \stc \; P_{sas'}^{0}>0,\\
    & \kappa_{s'} = 0, \forall \; s' \in \X \; \stc \; P_{sas'}^{0} = 0,\\
    & \bm{\kappa} \in \R^{\X}
\end{cases}\]
and the constraints on the right-hand side above are conic constraints. \jgc{Note that this equivalence remains true at $\xi=0$.}
Plugging this into~\eqref{eq:reformulation-concise-final}, we obtain the following corollary. 
\begin{corollary}\label{cor:reformulation-concise-kl}
Assume that the uncertainty set $\U$ is sa-rectangular and each set $\U_{sa}$ is of the form~\eqref{eq:uncertainty-kl}. Then the convex formulation~\eqref{eq:reformulation-concise-final} can be reformulated as follows:
\jgc{
\begin{equation}\label{eq:reformulation-concise-kl-conic}
    \begin{aligned}
    \max \quad  & \sum_{s \in \X} x_{s} \\
    \stc \quad & x_{s} \leq \sum_{a \in \A} \omega_{b}(s,a) \left(
   y_{sas'} + z_{sas'}  + \left(\frac{1-\gamma}{\gamma}\right)  u_{sa} + 
 \left(\frac{1+\log(\gamma)}{\gamma}\right)\alpha_{sa} - \xi_{sa}\theta_{sa} - q_{sa}\right), \forall \; s,s' \in \X,\\
 & (x_{s'},\alpha_{sa},z_{sas'}) \in \Kexp, \forall \; (s,a,s') \in \X \times \A \times \X, \\
 & (1,\alpha_{sa},u_{sa}) \in \Kexp, \forall \; (s,a) \in \X \times \A, \\
   & \sum_{s' \in \X}\kappa_{sas'} \leq \xi_{sa}, \forall \; (s,a) \in \X \times \A,\\
    & (\kappa_{sas'}/P^{0}_{sas'}, \xi_{sa},y_{sas'} - q_{sa})  \in \Kexp, \forall \; (s,a) \in \X \times \A, \forall \; s' \in \X \; \stc \; P^{0}_{sas'}>0, \\
    & \kappa_{sas'} = 0, \forall \; (s,a) \in \X \times \A, \forall \; s' \in \X \; \stc \; P^{0}_{sas'}=0, \\
    & x_{s} \geq 1, \forall \; s \in \X, \\
    &   \bm{\xi}_{sa}\in \R_{+}^{m},\alpha_{sa} \in \R_{+},  u_{sa} \in \R, \bm{z}_{sa} \in \R^{\X},\bm{y}_{sa} \in \R^{\X}, \bm{\kappa}_{sa} \in \R^{\X}, q_{sa} \in \R,  \bm{x} \in \R^{\X}, \forall \; (s,a) \in \X \times \A.
\end{aligned}
\end{equation}
}
\end{corollary}
}

\tb{
\section{Proofs for Section~\ref{sec:s-rectangular}}\label{app:proof-s-rectangular}
\proof{Proof of Lemma~\ref{lem:t-tilde-simplified-s-rec}}
Let $(s,a) \in \X \times \A$ and define $\hat{f}_{sa}\colon \R^{\X \times \A} \rightarrow \R$ such that $\hat{f}_{sa}(\bm{P}_{s}) = f(\bm{P}_{sa})$ for any $\bm{P}_{s} \in \R^{\X \times \A}$ and for $f$ defined in~\eqref{eq:definition-f-p}. With these notations, we have
\[ \tilde{t}(\bm{x})_{s} = \min_{\bm{P}_{s} \in \U_{s}} \sum_{a \in \A} \omega_{b}(sa)\hat{f}_{sa}(\bm{P}_{s}).\]
We provide a maximization formulation for $\tilde{t}(\bm{x})_{s}$. Using convex duality, we have that
\begin{align*}
  \min_{\bm{P}_{s} \in \U_{s}} \quad & \sum_{a \in \A} \omega_{b}(sa)\hat{f}_{sa}(\bm{P}_{s})
  \\
  = \quad &\\
  \max \quad& \sum_{a \in \A} \theta_{sa} + \min_{\bm{P} \in \R^{\X \times \A}} \left(\bm{\xi}_{s}\tr\bm{g}_{s}\right)(\bm{P}_{s}) +  \sum_{a \in \A} \omega_{b}(sa)\hat{f}_{sa}(\bm{P}_{s}) - \bm{P}_{sa}\tr(\bm{\mu_{sa}} + \theta_{sa} \bm{e}) \\
  \stc \quad  & \bm{\mu}_{sa} \in \R^{\X}_{+}, \theta_{sa} \in \R, \forall \; a \in \A, \\
    & \bm{\xi}_{s} \in \R^{m}_{+}
\end{align*}
This shows that 
\begin{align*}
  \min_{\bm{P}_{s} \in \U_{s}} \quad& \sum_{a \in \A} \omega_{b}(sa)\hat{f}_{sa}(\bm{P}_{s}) \\
  = \quad&\\
  \max \quad& \sum_{a \in \A} \theta_{sa} - \left( \left(\bm{\xi}_{s}\tr\bm{g}_{s}\right) + \sum_{a \in \A} \omega_{b}(sa)\hat{f}_{sa} \right)^{*}((\bm{\mu}_{sa} + \theta_{sa}\bm{e})_{a \in \A}) \\
   \stc \quad  & \bm{\mu}_{sa} \in \R^{\X}_{+}, \theta_{sa} \in \R, \forall \; a \in \A, \\
    & \bm{\xi}_{s} \in \R^{m}_{+}
\end{align*}
Using the same methods as in Section~\ref{sec:tractable counterparts - sa-rec} to eliminate the variables $\bm{\mu}_{sa}$ and $\theta_{a}$, and point (ix) from proposition 1.3.1 in \cite{hiriart1996convex}, we obtain that $\min_{\bm{P}_{s} \in \U_{s}} \sum_{a \in \A} \omega_{b}(sa)\hat{f}_{sa}(\bm{P}_{s})$ can be reformulated as
\begin{align*}
       \max \quad & - \left(\bm{\xi}_{s}\tr\bm{g}_{s}\right)^{*}(\bm{y}_{s})  +   \sum_{a \in \A} \omega_{b}(s,a) \left( \min_{s' \in \X} \{\alpha_{sa}\log(x_{s'}) + \frac{y_{sas'}}{\omega_{b}(s,a)} \} -\frac{\alpha_{sa}}{\gamma} \log \left(\frac{\alpha_{sa}}{\gamma}\right) +
 \frac{\alpha_{sa}}{\gamma}  \right) \\
    & \bm{\xi}_{s} \in \R^{m}, \bm{y}_{sa} \in \R^{\X}, \alpha_{sa} \in \R_{+}, \forall \; a \in \A.
    \end{align*}
\hfill \halmos
\endproof
}

\end{APPENDICES}

\end{document}